# Conforming, non-conforming and non-matching discretization couplings in discrete fracture network simulations

Alessio Fumagalli[a], Eirik Keilegavlen[a], Stefano Scialò[b,c]

[a]*Department of Mathematics, University of Bergen, Allègaten 41, Bergen 5007, Norway*
[b]*Dipartimento di Scienze Matematiche, Politecnico di Torino, Corso Duca degli Abruzzi 24, 10129 Torino, Italy*
[c]*Member of the INdAM research group GNCS*

**Abstract**

Simulations of fluid flow in naturally fractured rocks have implications for several subsurface applications, including energy storage and extraction, and waste storage. We are interested in flow in discrete fracture networks, which explicitly represent flow in fracture surfaces, but ignore the impact of the surrounding host rock. Fracture networks, generated from observations or stochastic simulations, will contain intersections of arbitrary length, and intersection lines can further cross, forming a highly complex geometry. As the flow exchange between fractures, thus in the network, takes place in these intersections, an adequate representation of the geometry is critical for simulation accuracy. In practice, the intersection dynamics must be handled by a combination of the simulation grid, which may or may not resolve the intersection lines, and the numerical methods applied on the grid. In this work, we review different classes of numerical approaches proposed in recent years, covering both methods that conform to the grid, and non-matching cases. Specific methods considered herein include finite element, mixed and virtual finite elements and control volume methods. We expose our methods to an extensive set of test cases, ranging from artificial geometries designed to test difficult configurations, to a network extruded from a real fracture outcrop. The main outcome is guidances for choice of simulation models and numerical discretization with a trade off on the computational cost and solution accuracy.

domain decomposition

## 1. Introduction

Flow through fractured rocks is important for several applications including $CO_2$ storage [54, 49], geothermal energy recovery [55], and nuclear waste disposal [56]. Several classes of mathematical models have been developed to represent flow through fractured rocks. These can be characterized according to the degree in which they represent fractures explicitly, or by substitution by average flow properties [12, 23]. Here, we focus on one extreme case, the discrete fracture network (DFN) approach, where any flow in fractures not represented in the simulation model is ignored. Such methods are mainly applied when the host rock has negligible permeability outside fractures, for instance in low-permeable sandstone, shales, or granite. DFN models have successfully been applied, for instance, in modeling of $CO_2$ storage [40] and hydraulic fracturing [45].

Natural fracture networks can exhibit highly complex geometric configurations even in two-dimensional outcrops [36, 29], and the inclusion of the third dimension further increases complexity [24]. Due to the lack of access, hard data on subsurface fractures is hard to come by, and fracture networks for simulation models are therefore commonly created by stochastic generation, informed by bore hole data, seismic imaging, outcrop analogues and other sources, see *e.g.* [42, 50, 26].

The generation of stochastic fracture networks is challenging by itself [25]. Possible approaches range from the drawing of sets of fractures from statistical distributions, with no conditioning from the already generated fractures, to more advanced ones that attempt to honor relation between different geological objects [21, 58, 43, 52]. As stochastic generation does not necessarily aim to mimic the actual fracturing process, the generated networks may not be realistic in a geological sense, as measured by geometric configurations, such as intersection types, and angles and distances between intersection lines. Nevertheless, the generated networks are believed to give realistic representation of larger-scale flow properties of the network, such as permeability and breakthrough time, and flow simulation models should be designed to handle general network configurations, based on hard data or stochastic realizations.

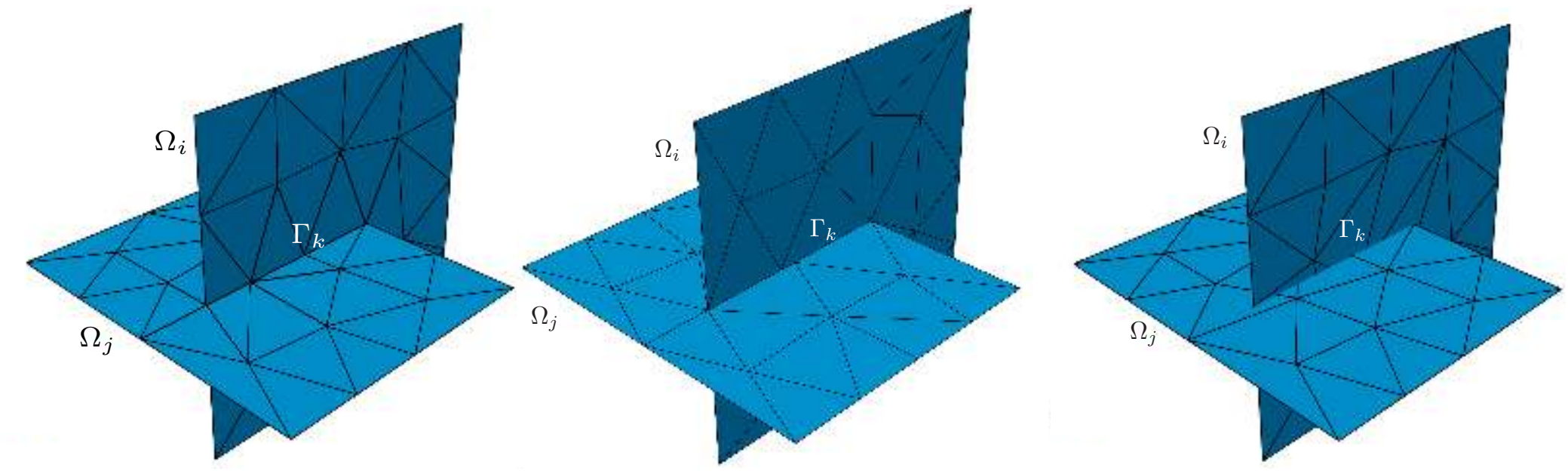

Figure 1: Conforming (left), triangular non conforming (center) and non matching mesh (right) for two intersecting fractures

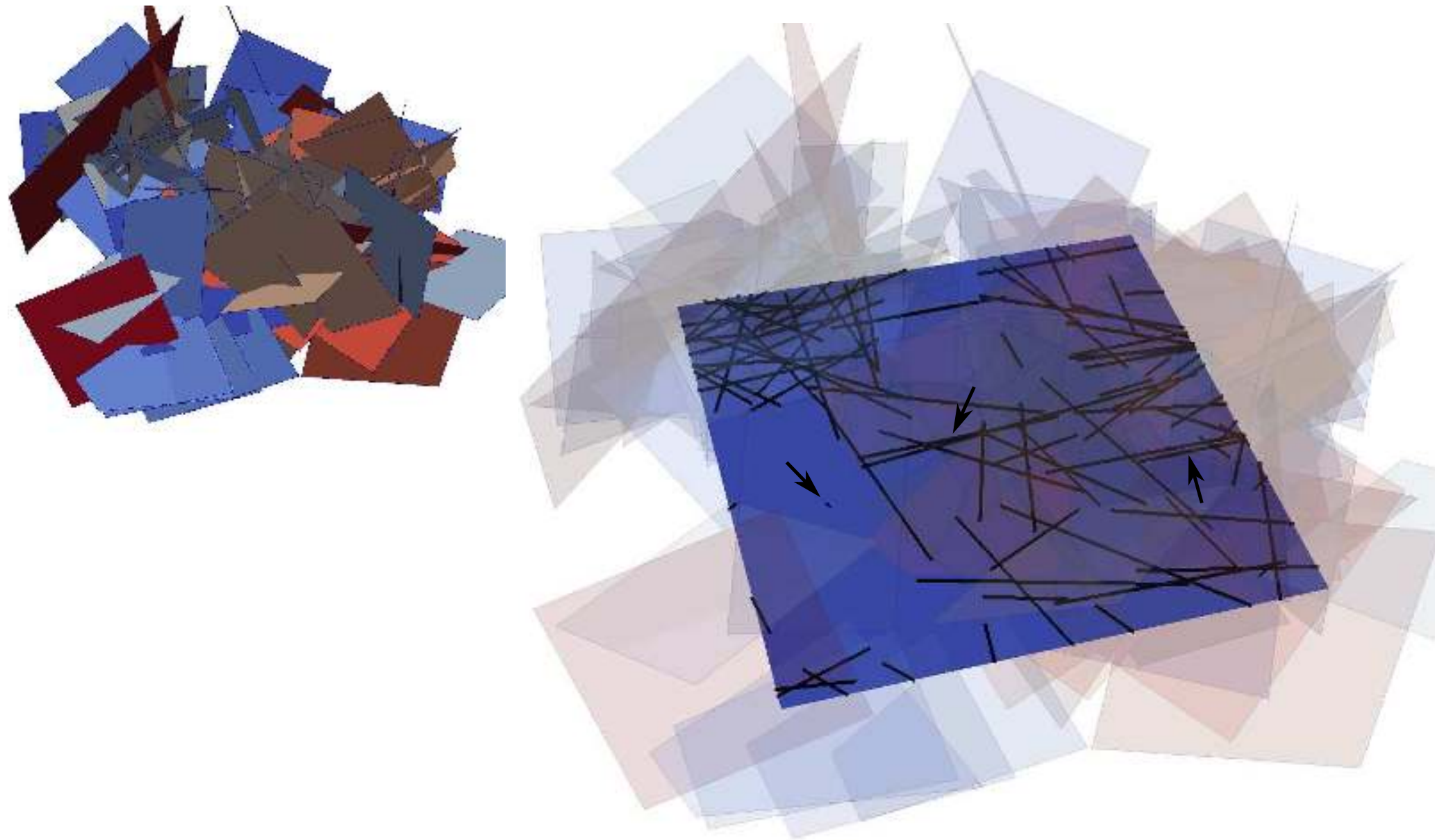

Figure 2: Example of hard-to-mesh geometrical features: traces forming small angles, parallel close traces, traces of different length scales

The geometric complexity of the fracture network poses the main challenge to numerical simulation models based on DFN principles. As flow can only travel between fractures via fracture intersections, it is paramount that flow through the intersections is represented in the discretized model. This puts practical constraints on either the computational mesh, or the numerical method, or both, and gives rise to what can be considered two classes of numerical methods for DFN flow models. The first class requires the computational mesh to conform to all intersection lines, with no hanging nodes, see Figure 1 (left) for an example. With the intersection lines explicitly represented as internal boundaries in the fracture surfaces, discretization of the flow equations takes the form of a set of standard flow discretizations, coupled by boundary conditions. Depending on which model is chosen for dynamics at the interface, this discretization can be made by most standard numerical methods that can handle pressure and flow boundary conditions, examples relevant in this setting include standard finite volume and finite element approaches. The fully conforming approach moves all complexity to the meshing algorithm, with each of the configurations in Figure 2 posing challenges that can only partly be resolved by mesh refinement. To ensure node conformity from all fracture surfaces, the fractures must be meshed simultaneously, see [38, 39, 46, 47] for some approaches. The difficulties in terms of meshing can partly be alleviated by meshing the fracture surfaces independently, without requiring matching nodes along intersections, see Figure 1 (center). This partly conforming approach reduces the meshing to a standard two-dimensional problem for which high quality software is readily available, *e.g.* [34, 53]. However, the choice of numerical methods is limited to approaches that can handle hanging nodes, of particular interest to us here is the Virtual Element Method, (VEM) [3, 4, 5, 8, 13].

The second class of numerical schemes removes the requirement of mesh conformity to intersection lines completely, see Figure 1 (right) for an example. The meshing then reduces to the relatively easy task of griding a decoupled set of fracture surfaces. However, the matching of pressure and fluxes in intersecting fractures is now left to the discretization scheme. As the intersection lines arbitrarily cut cells, fully or partly, specialized numerical methods are needed. Possible options include finite elements with the intersection treated as internal jumps and non-conformal discretizations [22, 30], with mass conservation imposed by an optimization scheme [14], and using virtual elements, [10].

In this paper, we offer a comprehensive computational comparison of nu-merical discretization schemes for DFN models. Our goal is to study the trade off between efficient meshing and the use of standard numerical methods, and to provide guidance for the choice of simulation models. To that end, we consider in total 7 numerical methods with different order of polynomial degrees, including members tailored for fully conforming, partly conforming and non-conforming meshes. The test cases include setups designed to probe performance on what is known to be difficult geometric configurations. The geometries and the computational mesh are available to the reader to be used for further investigations. Some of the results are produced by the PorePy library, see [41] and `http://github.com/pmgbergen/porepy` for more details, others using the Geoscore++ library, see `https://areeweb.polito.it/geoscore/software`.

The paper is organized as follows. The mathematical models to describe the fluid flow in the DFN are presented in Section 2, in both primal and dual formulation. In Section 3 we review the possible discretization strategies to tackle DFNs, in particular for conforming, non-conforming, and non-matching methods. Extensive experimental results are presented in Section 4 to validate in different configurations the considered numerical schemes. Conclusions follow in Section 5.

## 2. Mathematical model

In this section we introduce the mathematical model able to describe a single-phase flow in DFNs. We consider the fractures as (approximated by) planar objects, following *e.g.* [37, 20, 32, 2, 27]. The extension to curved fractures is possible, but the analysis and the presentation will increase in complexity, see [28, 35].

Let us consider a convex flat domain $\Omega_i \subset \mathbb{R}^3$ embedded in the three-dimensional space with boundary $\partial\Omega_i$. $\Omega_i$ represents one fracture with index $i \in \{1, \dots, N_\Omega\}$, where $N_\Omega$ is the total number of fractures. We suppose that the fractures are non-overlapping, *i.e.* $\lambda_2(\cap_i \Omega_i) = 0$ for all $i \in \{1, \dots, N_\Omega\}$ with $\lambda_2$ the bi-dimensional Lebesgue measure. We indicate by $\Omega = \cup_i \Omega_i$ for all $i \in \{1, \dots, N_\Omega\}$ the DFN with boundary $\partial\Omega$. Given $\Omega_i$ and $\Omega_j$ we indicate a trace (a fracture intersection) the line $\Gamma_k = \Omega_i \cap \Omega_j \neq \emptyset$. A natural order of indexes, related to the fractures involved, can be set for the traces $k \in \{1, \dots, N_\Gamma\}$ by a function $t : \mathbb{N} \times \mathbb{N} \to \mathbb{N}$ such that $k = t(i, j)$ with $\Gamma_k = \Omega_i \cap \Omega_j$. Note that $t(i, j) = t(j, i)$, also $N_\Gamma$ indicates the total number of traces. We suppose that a trace is formed only by the intersection of two fractures, however the extension to multiple fracture intersection is immediate. Finally, we indicate by $\Gamma = \cap_k \Gamma_k$ for all $k \in \{1, \dots, N_\Gamma\}$ the set of all traces and by $\Gamma_{\Omega_i}$ the set of traces on fracture $\Omega_i$, $i \in \{1, \dots, N_\Omega\}$. Figure 3 represents a simple DFN with the introduced notations.

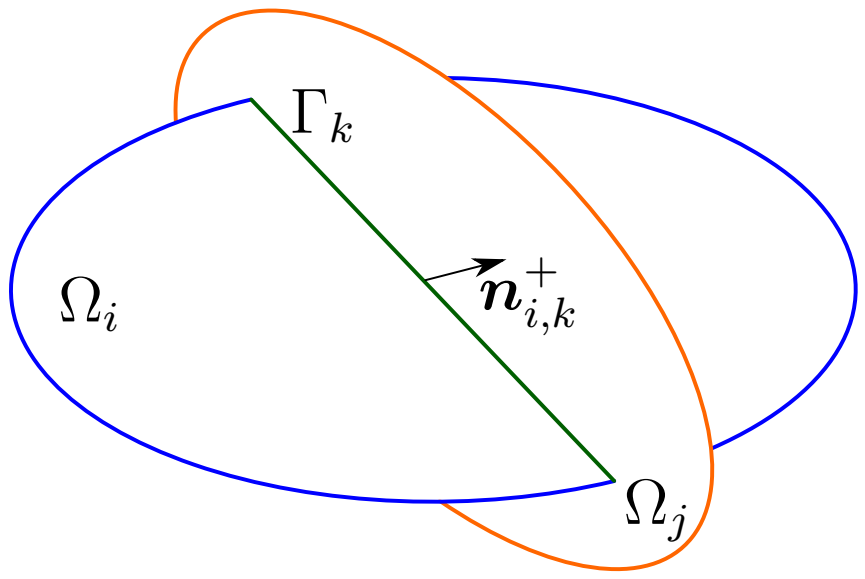

Figure 3: Graphical representation of a DFN. In the picture we have assumed that $k = t(i,j)$.

In this work we are interested to model an incompressible single-phase flow, we consider thus a generalization of the classical Darcy problem for the description of the hydraulic head $h$ and the flux $\boldsymbol{u}$ in the DFN. We indicate by a subscript in $\{1, \dots, N_\Omega\}$ the restriction of $h$ and $\boldsymbol{u}$ on a single fracture. On a single fracture the model reads

$$\begin{aligned} \boldsymbol{u}_i &= -K_i \nabla h_i \\ \nabla \cdot \boldsymbol{u}_i &= f_i \end{aligned} \quad \text{in } \Omega_i \setminus \Gamma. \tag{1a}$$

Here $K_i \in \mathbb{R}^{2\times 2}$ is the symmetric and positive defined permeability matrix on the tangent plane of $\Omega_i$ and $f_i$ is a scalar source/sink term. The differential operators considered in (1a), and the flux field $\boldsymbol{u}$ are defined on the tangent space of $\Omega_i$. The considered model does not explicitly include the fracture aperture as parameter however, following the idea presented in [44, 30, 51, 32], fracture aperture can be included in $K_i$ by scaling.

To allow for discontinuities over traces within a fracture, and thus flow between fractures, the values of the hydraulic head and the normal component of the flux are doubled at $\Gamma_k$. Let $\boldsymbol{n}^*_{i,k}$, represent the outward unit normal at $\Gamma_k$ tangent to $\Omega_i$ for each side $* \in \{+, -\}$ (left and right). At each trace we consider the continuity of the normal flux and the hydraulic head, the following coupling conditions hold

$$\begin{aligned} \sum_{*\in\{+,-\}} \boldsymbol{u}_i \cdot \boldsymbol{n}^*_{i,k}|_{\Gamma^*_k} + \boldsymbol{u}_j \cdot \boldsymbol{n}^*_{j,k}|_{\Gamma^*_k} = 0 & \\ & \quad \text{on } \Gamma_k \\ h_i|_{\Gamma^+_k} = h_j|_{\Gamma^+_k} = h_i|_{\Gamma^-_k} = h_j|_{\Gamma^-_k} & \end{aligned} \tag{1b}$$

with $(i,j) = t^{-1}(k)$. The first condition in (1b) represents conservation of mass, while the second condition enforces continuity of the hydraulic head. In (1b) we consider also $\cdot|_{\Gamma^*_k}$ a suitable trace operator which, with an abuse of notation, we have used for both the flux and the hydraulic head for each side $* \in \{+,-\}$. In cases where the trace coincides with the boundary of a fracture, as happens in an T-, L- or Y-type intersection, (1b) is modified to only account for one side of the fracture. The coupling conditions (1b) can be viewed as an approximation of the general conditions described in [32, 19, 33] when the aperture of the intersection is of negligible size. The considered hypothesis is reasonable for many applications.

For each fracture and the full DFN, we split the boundary into two non-overlapping portions named $\Gamma_{i,D}$ and $\Gamma_{i,N}$, for the fracture $\Omega_i$, and $\Gamma_D$ and $\Gamma_N$ for the DFN $\Omega$. On $\Gamma_{i,D}$ and $\Gamma_D$ we prescribe boundary conditions for the hydraulic head and on $\Gamma_{i,N}$ and $\Gamma_N$ for the flux, namely, for a single fracture:

$$\begin{aligned} h_i|_{\Gamma_{i,D}} &= \overline{h}_i && \text{on } \Gamma_{i,D} \\ \boldsymbol{u}_i \cdot \boldsymbol{n}_{iN}|_{\Gamma_{i,N}} &= \overline{u}_i && \text{on } \Gamma_{i,N} \end{aligned} \tag{1c}$$

where $\cdot|_{\Gamma_{i,D}}$ and $\cdot|_{\Gamma_{i,N}}$ denote suitable trace operators, $\boldsymbol{n}_{i,N}$ is the unit normal of $\partial\Omega_i$ pointing outward with respect to $\Omega$, and $\overline{h}_i$ and $\overline{u}_i$ are given data.

The problem in mixed form is: find $(h, \boldsymbol{u})$ such that the system described in (1) holds true. The problem can be re-written in primal form, with only the hydraulic pressure as unknown, combining the two equations in (1a) into

$$-\nabla \cdot K_i \nabla h_i = f_i \quad \text{in } \Omega_i \setminus \Gamma, \tag{2a}$$

with coupling conditions, derived from (1b), between intersecting fractures

$$\begin{aligned} \sum_{*\in\{+,-\}} K_i \nabla h_i \cdot \boldsymbol{n}^*_{i,k}|_{\Gamma^*_k} + K_i \nabla h_j \cdot \boldsymbol{n}^*_{j,k}|_{\Gamma^*_k} = 0 & \\ & \quad \text{on } \Gamma_k. \\ h_i|_{\Gamma^+_k} = h_j|_{\Gamma^+_k} = h_i|_{\Gamma^-_k} = h_j|_{\Gamma^-_k} & \end{aligned} \tag{2b}$$

Finally, the boundary conditions expressed in (1c) can be recast as

$$\begin{aligned} h_i|_{\Gamma_{i,D}} &= \overline{h}_i && \text{on } \Gamma_{i,D} \\ -K_i \nabla h_i \cdot \boldsymbol{n}_{i,N}|_{\Gamma_{i,N}} &= \overline{u}_i && \text{on } \Gamma_{i,N} \end{aligned}. \tag{2c}$$

The problem in primal form is: find $h$ such that the system described in (2) holds true. Under reasonable hypotheses on the data, the weak formulations of the primal and mixed formulations are well posed, *e.g.* [30, 32, 51, 11].

## 3. Numerical discretization

In this section we recall three common strategies to derive a numerical scheme for the solution of the flow in DFNs. The main concern is how the matching conditions (1b) and (2b) along the traces are implemented, different choices will lead to different advantages and drawbacks. The principal difference is how the connections are handled at fracture intersections: conforming, non-conforming, and non-matching. These are discussed in Subsections 3.1, 3.2, and 3.3, respectively. In each subsection appropriate numerical methods are briefly discussed; a detailed description of the methods can be found in the references.

### *3.1. Conforming discretization at traces*

In this part we introduce the key concepts for conforming discretizations for problems (1) and (2). Conforming discretizations exploit the fact that the considered meshes are conforming at each trace in the network. Figure 1 on the left shows two discretized fractures and their meshes. These are conforming in the sense that the nodes at the trace form a contiguous set of edges which are part of elements in the two fracture meshes. The figure shows simplices but other geometries are allowed.

For a chosen discretization, the implementation of the coupling conditions (1b) or (2b) is normally done by considering, if present, multiple degrees of freedom associated with the edges or points at the traces and the Lagrange multipliers technique to enforce the coupling conditions. The procedure considers the traces as part of the boundary and new equations couple the associated degrees of freedom. The actual implementation of this procedure is specific to the numerical scheme, we give two examples related to our case.

If we consider a numerical scheme that solves directly (1), the continuity of the normal fluxes at each trace is enforced by Lagrange multipliers while the continuity of the hydraulic head is weakly given by integration by part. The system of equation in weak form becomes: find $(\boldsymbol{u}_i, h_i, \lambda_i) \in V_i \times Q_i \times \Lambda_k$ in each $\Omega_i$ and $\Gamma_k$ such that

$$\begin{aligned}
&\left(K_i^{-1}\boldsymbol{u}_i, \boldsymbol{v}_i\right)_{\Omega_i} - (h_i, \nabla \cdot \boldsymbol{v}_i)_{\Omega_i} + \sum \left\langle \lambda_k, \boldsymbol{v}_i \cdot \boldsymbol{n}_{i,k}^*|_{\Gamma_k^*} \right\rangle = 0, && \forall \boldsymbol{v}_i \in V_i;\\
&- (\nabla \cdot \boldsymbol{u}_i, q_i)_{\Omega_i} = - (f_i, q_i)_{\Omega_i} && \forall q_i \in Q_i\\
&\sum \left\langle \boldsymbol{u}_i \cdot \boldsymbol{n}_{i,k}^*|_{\Gamma_k^*}, \mu_k \right\rangle = 0, && \forall \mu_k \in \Lambda_k
\end{aligned}$$

the summations are done on $k \in \Gamma_{\Omega_i}$ and $* = \{+, -\}$, the latter indicating each side of the $k$-trace on the $i$-fracture. The Hilbert spaces $V_i$, $Q_i$, and $\Lambda_k$ are chosen accordingly to ensure well-posedness of the problem. The duality pairing can be identified by $L^2$ scalar products, since a higher regularity ($L^2(\Gamma_k)$ instead of $H^{-\frac{1}{2}}(\Gamma_k)$) is normally assumed for the trace of the normal fluxes, see, *e.g.*, [32] for a more detailed discussion on this aspect. The method MVEM-CONF (mixed-VEM for conforming mesh), presented in the aforementioned work, fits in this framework and approximate elements of $Q_i$ as piecewise constant on each elements, $V_i$ as a suitable VEM space, and $\Lambda_k$ as piecewise constant on each edge of the traces. Also, the simplicial meshes considered for MVEM-CONF can be coarsened by gluing neighbor cells of small volume, as discussed in [33]. We refer to this approach as MVEM-COARSE.

To approximate (2) in the case of cell-centered finite volume methods, the strategy is to reconstruct the fluxes at the traces and enforce their conservation. Also in this case the continuity of the hydraulic head is weakly given by integration by part. We consider a tessellation $\mathcal{T}$ on the fractures. For each cell $c$ its set of edges is $\mathcal{E}(c) \subset \mathcal{E}$, with $\mathcal{E}$ the set of edges of $\mathcal{T}$, and for each edge $e$ its set of neighboring cells is $\mathcal{T}(e)$. If an edge belongs to a trace, the neighboring cells are part of different fractures. Problem (2) becomes

$$\begin{aligned}
&\sum_{e \in \mathcal{E}(c)} F_{c,e}(h) = I(f_i|_c) && \forall c \in \mathcal{T}\\
&\sum_{c \in \mathcal{T}(e)} F_{c,e}(h) = 0 && \forall e \in \mathcal{E}
\end{aligned}$$

where $I(f_i|_c)$ is a numerical approximation of the source term and $F_{c,e}$ represents the outflow from the cell $c$ through the edge $e$. Each numerical scheme approximates the flux $F_{c,e}$ by the values of the hydraulic head on a suitable patch $\sigma(F_{c,e})$ of cells. In our application two schemes are considered: TPFA

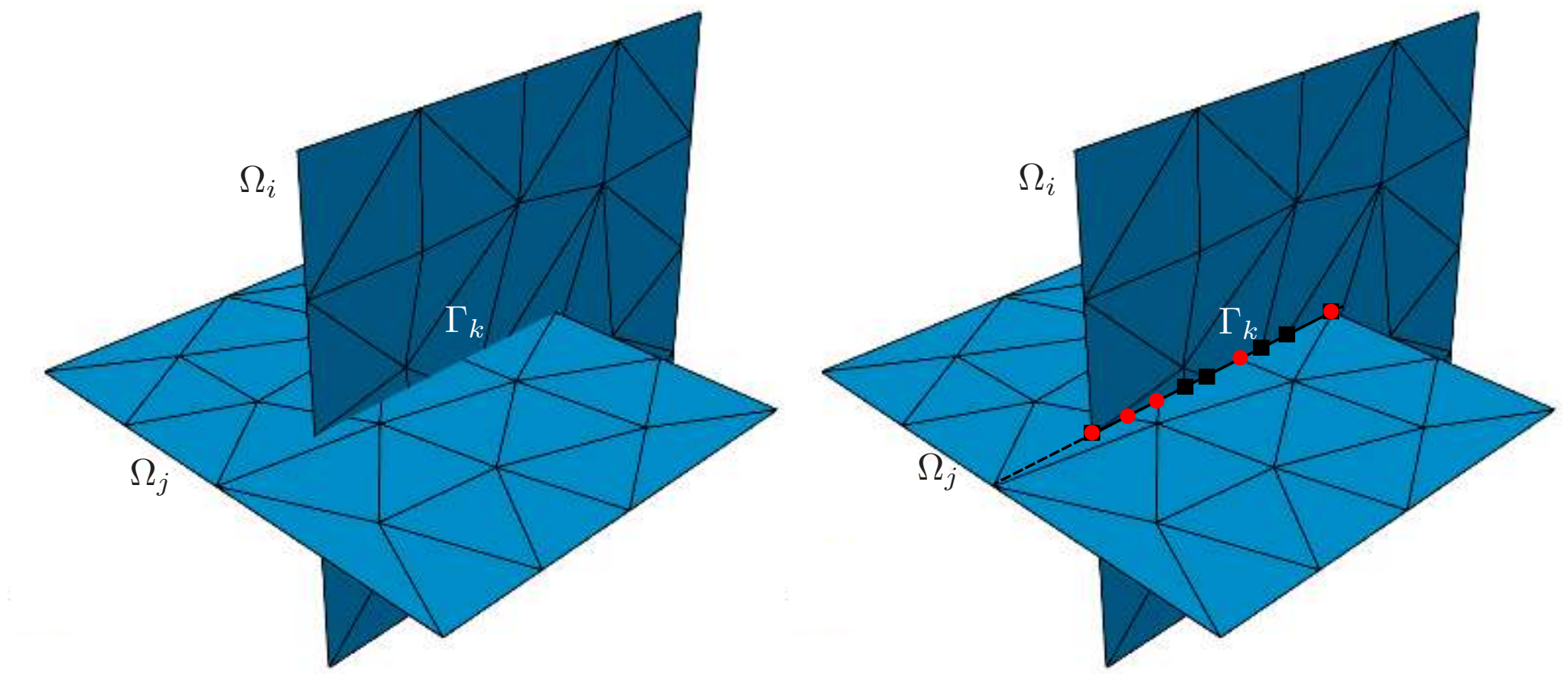

Figure 4: Non matching mesh (left) and corresponding polygonal non-conforming mesh (right). Mesh elements are cut into polygons conforming to the traces of each fracture. Black squares mark mesh element vertexes on $\Gamma_k$ for the mesh on $\Omega_i$, red dots for $\Omega_j$

(two-point flux approximation) where the patch $\sigma(F_{c,e})$ is edge based, and MPFA (multi-point flux approximation) where the patch $\sigma(F_{c,e})$ is node based [1]. The former suffers from inconsistency on non-$K$ orthogonal tessellations, which will be the case for our grids. This is remedied by MPFA, although our implementation still applies a two-point approximation for flow into traces.

### *3.2. Non-conforming discretization at traces*

To overcome some of the difficulties to generate conforming triangular meshes on intricate network, it is possible to generate polygonal conforming meshes which require less constraints in their construction. Let us start from triangular meshes, independently built on each fracture in $\Omega$, and thus possibly non conforming to the traces, as the one shown in Figure 4 (left). Then the triangles of the mesh are cut, on each fracture, into sub-polygons conforming to the traces, as shown in Figure 4 (right). This is a conforming mesh of each fracture individually, but it is not a conforming mesh of the whole network, as the mesh element vertexes of fracture $\Omega_i$ and $\Omega_j$ do not, in general, match on trace $\Gamma_k$, with $(i,j) = t^{-1}(k)$, as shown in Figure 4 (right), where black squares mark the vertexes of the mesh on $\Omega_i$ on $\Gamma_k$ and red dots mark, instead the vertexes of the mesh on $\Omega_j$ on $\Gamma_k$. This non-conforming polygonal mesh can be used to solve problem (2), using the VEM to handle polygonal elements in conjunction with the mortar method to take care of the non conformity at the traces.

Let us introduce, for each fracture $\Omega_i$, $i \in \{1,\ldots,N_\Omega\}$, and for each trace $\Gamma_k$ in $\Gamma_{\Omega_i}$ the function $[u]_{i,k} = \sum_{*\in\{+,-\}} \boldsymbol{u}_i \cdot \boldsymbol{n}^*_{i,k}|_{\Gamma_k^*}$, corresponding to the jump of the co-normal derivative of $h_i$ across $\Gamma_k$ on $\Omega_i$. We have $[u]_{i,k} \in U_k$, with $U_k$ a suitable Hilbert space defined on $\Gamma_k$. Then, on each fracture $\Omega_i$, $i \in \{1,\ldots,N_\Omega\}$, functions $[u]_i \in U_{\Omega_i} = \prod_{k\in\Gamma_{\Omega_i}} U_k$ are defined as: $[u]_i = \prod_{k\in\Gamma_{\Omega_i}} [u]_{i,k}$. Let us choose the function space $V_i = H^1_0(\Omega_i)$, and thus let us set $U_k = H^{-\frac{1}{2}}(\Gamma_k)$. Assuming, for simplicity, homogeneous Dirichlet boundary conditions on all fracture boundaries, problem (2) can be written in weak formulation as: for all $i \in \{1,\ldots,N_\Omega\}$, find $h_i \in V_i$ and $[u]_i \in U_{\Omega_i}$ such that

$$\begin{aligned}
&(K_i\nabla h_i, \nabla v_i)_{\Omega_i} - \sum_{k\in\Gamma_{\Omega_i}} \sigma_i(k)\,\langle [u]_{i,k}, v_i|_{\Gamma_k}\rangle_{U_k,U_k'} = (f_i, v_i)_{\Omega_i}\,, && \forall v_i \in V_i;\\
&-\sum_{k\in\Gamma_{\Omega_i}} \sigma_i(k)\,\langle v_i|_{\Gamma_k}, \psi|_{\Gamma_k}\rangle_{U_k',U_k} = 0, && \forall \psi \in U_{\Omega_i}
\end{aligned} \tag{3}$$

being $\sigma_i(k) : [1, N_\Gamma] \mapsto \{-1, 1\}$ a function defined as follows:

$$\sigma_i(k) = \begin{cases} -1 & \text{if } i = \min\left(t^{-1}(k)\right) \\ 1 & \text{otherwise} \end{cases}$$

thus selecting the master and slave fracture for each couple of fractures meeting at each trace. System (3) is common to other methods that are designed to solve problem (2) by using the mortar technique [6, 57].

Once system (3) is discretized, it can be re-written in matrix form giving a classical saddle-point algebraic system that is solved to obtain the discrete solution. For our study, the resulting numerical schemes will be termed VEM-M (mortar VEM in primal formulation) of order $n = 1, 2$. Details can be found in [7, 9].

The flexibility of the VEM in handling polygonal elements also with flat angles at some vertexes can be exploited in order to transform the non-conforming mesh of Figure 4 (right) into a conforming mesh of the whole network. This can be easily done by adding the vertexes of the mesh elements of fracture $\Omega_i$ that lie on trace $\Gamma_k$ to the elements of fracture $\Omega_j$ on $\Gamma_k$, and vice versa, with $k = t(i, j)$ for each trace $k \in \{1,\ldots,N_\Gamma\}$. With reference to Figure 4 (right), this means that now the mesh elements of fracture $\Omega_i$ have both the vertexes marked with black squares and red dots, and similarly for $\Omega_j$. The resulting numerical scheme retains the saddle-point structure of problem (3), but it is no-more necessary to resort to the mortar method to enforce the continuity condition, which can be now imposed simply equating the matching degrees of freedom of the fractures meeting at each trace by means of Lagrange multipliers. The resulting schemes are termed VEM-C (VEM in primal formulation with conforming mesh) of order $n = 1, 2$. Details can be found in [11, 9]. This approach can be seen as conforming method where the potentialities of the numerical scheme are exploited at best.

The same approaches can be used to solve problem (1), using the mixed VEM. This gives the numerical schemes denoted as MVEM of order $n = 0, 1, 2$. Details can be found in [11]. The main difference between MVEM and MVEM-CONF is the way the grids are constructed. This will impact on the numerical results as presented in Section 4.

### *3.3. Non-matching discretization at traces*

The generation of a mesh that matches fracture intersection might not be trivial for intricate networks with a large number of fractures. Moreover, randomly generated networks might display geometrical features, such as narrow angles between intersecting traces, very small traces, or non intersecting traces on a fracture running parallel and very close to each other, and the generation of a mesh conforming to such features would require a very large number of mesh elements, independently of the required mesh-size. In order to overcome such difficulties, a numerical scheme is proposed in [14, 15, 18] for problem (2) that relaxes any conformity requirement of the mesh at the traces. The method relies on a cost functional to express the error in the fulfillment of the matching conditions (2b), and the solution found as the minimum of this functional, constrained to satisfy the Darcy law on the fractures. The method is here recalled in its simplest form, also assuming homogeneous Dirichlet boundary conditions on the whole boundary of a DFN $\Omega$. Let us then introduce the function spaces, defined as previously as $V_i := \mathrm{H}^1_0(\Omega_i)$, $i \in \{1, \dots, N_\Omega\}$ and $U_k := \mathrm{H}^{-\frac{1}{2}}(\Gamma_k)$, for all $k \in \{1, \dots, N_\Gamma\}$. The cost functional $J$ is defined as follows:

$$J(h, [u]) = \frac{1}{2} \sum_{k=1}^{N_\Gamma} \|h_i - h_j\|^2_{U'_k} + \|[u]_{i,k} + [u]_{j,k}\|^2_{U_k}$$

with $[u] = \prod_{i=1}^{N_\Omega}[u]_i$. Problem (2a) can be written in weak formulation on each fracture as: for all for $i \in \{1,\ldots,N_\Omega\}$, find $h_i \in V_i$ such that

$$(K_i \nabla h_i, \nabla v_i)_{\Omega_i} = (f_i, v_i)_{\Omega_i} + \sum_{k \in \Gamma_{\Omega_i}} \langle [u]_{i,k}, v_i \rangle_{U_k, U'_k}, \quad \forall v_i \in V_i. \tag{4}$$

As shown in the above references, it is proven that the solution of problem (2) is equivalent to the following optimization problem:

$$\min J(h, [u]) \text{ subject to (4).} \tag{5}$$

Introducing a triangulation $\mathcal{T}_{\delta,i}$ independently on each fracture $\Omega_i$, with $\delta$ the mesh parameter, we define a finite element space $V_{\delta,i} \subset \mathrm{V}_i$ associated to $\mathcal{T}_{\delta,i}$. Similarly, for $k \in \{1,\ldots,N_\Gamma\}$, $\mathcal{T}_{\delta,i,\Gamma_k}$ is the triangulation on $\Gamma_k$ related to fracture $\Omega_i$, which is in general different from the discretization of $\Gamma_k$ on $\Omega_j$, $(i,j) = t^{-1}(k)$, and $U_{\delta,i,k} \subset U_{i,k}$ is the associated finite element space. Then, it is possible to rewrite the minimization problem (5) in the sub-spaces $V_{\delta,i}$, $U_{\delta,i,k}$, $i \in \{1,\ldots,N_\Omega\}$, $k \in \{1,\ldots,N_\Gamma\}$, and collect the integrals of basis functions into coefficient matrices. This produces a saddle-point problem, whose Shur-complement matrix is positive definite [10, 15], and can be solved via a gradient based approach.

The simplest discretization of this algorithm, using standard linear Lagrangian finite elements, will be referred to as OPT-FEM. However, the solution displays jumps of the co-normal derivative across the traces corresponding to the flux entering or leaving the fracture through the trace. To better approximate this, the discrete space $V_{\delta,i}$ on each fracture $\Omega_i \in \Omega$ can be enriched with suitable basis function, following the XFEM [31]. We will refer to this second approach as OPT-XFEM.

### *3.4. Computational frameworks*

To summarize, we consider 7 numerical schemes, belonging to one of the previous class of methods, for comparison in next Section 4, giving different complex configurations. To this purpose two software frameworks are developed. In particular TPFA, MPFA, MVEM-CONF, and MVEM-COARSE are developed using PorePy library [41]. For this reason the meshes considered by TPFA, MPFA, and MVEM-CONF in each simulation are the same and MVEM-COARSE applies a coarsening algorithm on them. Methods OPT-FEM,OPT-XFEM are implemented in Matlab® and the triangular non matching mesh is generated on each fracture independently using the software Triangle [53].

A parallel C++ version based on MPI and on CUDA [48] is available for the OPT-FEM method, [16, 17]. The same triangular mesh is then post processed to obtain a polygonal non-conforming mesh, used for methods VEM-M, or a polygonal conforming mesh for methods VEM-C and MVEM, as discussed above.

The labels used to identify the various methods are summarized in Table 1, where the following nomenclature is used: the formulation is either mixed (M) or primal (P), referring if the problem is considered with the formulation (1) or (2), respectively. Four different kind of meshes are proposed: triangular conforming to fracture intersections (T-C), triangular non-matching to intersections (T-NM), polygonal conforming (P-C), polygonal non-conforming (P-NC) and polygonal conforming deriving from a coarsening process, starting from a triangular conforming mesh ($P^\star$). Finally, specific strategy adopted to couple different fracture meshes are: functional based (Fb), mortar method (Mo), Lagrange multiplier (LM).

| *Label* | *Method* | *Form.* | *Mesh* | *Matching* |
|---|---|---|---|---|
| OPT-XFEM | Optimization & XFEM | P | T-NM | Fb |
| OPT-FEM | Optimization & FEM | P | T-NM | Fb |
| VEM-C# | VEM of order # | P | P-C | Mo |
| VEM-M# | VEM of order # | P-NC | P | Mo |
| MVEM# | MVEM of order # | M | P-C | LM |
| MVEM-CONF | MVEM order 0 | M | T-C | LM |
| MVEM-COARSE | MVEM of order 0 | M | $P^\star$ | LM |
| MPFA | Multi-point flux approx. | P | T-C | LM |
| TPFA | Two-point flux approx. | P | T-C | LM |

Table 1: Labels used to denote the various considered methods. *Formulation* (*Form.*): P = primal, M = mixed; *Mesh*: T-NM = triangular non-matching; T-C = triangular conforming; P-C = polygonal conforming; P-NC = polygonal non-conforming; $P^\star$ = polygonal from coarsening; *Matching*: Fb = functional based, Mo = mortar, LM = Lagrange multiplier.

## 4. Examples

In this section, we test the different numerical approaches on a set of test cases. In Subsection 4.1, we consider a convergence test on a simple case with three intersecting fractures with an explicit analytical solution. This allows us to benchmark the methods in a setting without geometric complexities. In Subsection 4.2 and 4.3 we consider three cases that are designed to test the methods for different geometric configurations, namely short traces and small angles between traces. These are geometries that are likely to arise in cases with many fractures, and the main motivation for the tests is to study method accuracy in an isolated setting as the level of geometric complexity is increased. For these cases, we do not carry out full convergence tests, but rather focus on accuracy, first on fine grids, establishing the potential performance each method, and second on relatively coarse grid that resemble what may be affordable in larger scale simulations. Accuracy is measured in terms of distribution of hydraulic head along selected lines, and also on the bulk flow properties of the network, e.g. the upscaled permeability. Finally, in Subsection 4.4 we consider a relatively complex network with 66 fractures, extruded from a real interpreted outcrop from Western Norway. Again we focus on method accuracy as the grid resolution is decreased from what can be termed highly resolved to much coarser cells.

A comparison of accuracy should take into account the computational cost of the numerical methods. The methods considered herein differ significantly in the formulations, ranging from standard cell-centered methods, via higher order and dual formulations to optimization based approaches. For simplicity, we use cell number as a proxy for computational cost, but it should be kept in mind that the number of unknowns, as well as the availability of efficient solvers and preconditioners will vary significantly between the methods.

In some of the following example we consider the Darcy flux for some comparisons. It is important to note that almost all the numerical schemes compute directly the fluxes, with the only exception of the method based on VEM-C, in which case it is necessary to post-process the hydraulic head to reconstruct the fluxes. In this latter case the reconstructed flux may suffer from lower accuracy.

### *4.1. Numerical evidence of error decay*

The first test is targeted at showing the convergence properties of the various proposed methods. Let us consider the domain $\Omega = \bigcup_{i=1}^{3} \Omega_i$ as shown in Figure 5, composed of three fractures with three fracture intersections. The following hydraulic head function is defined in $\Omega$, where a reference system

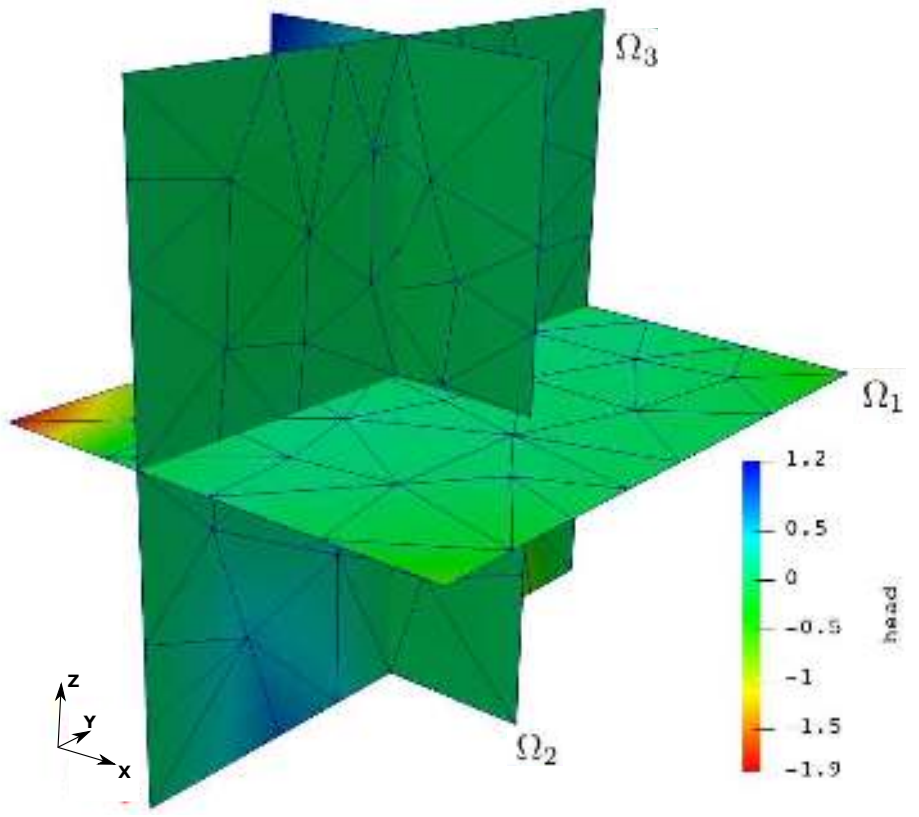

Figure 5: Computational domain and solution for the test problem of Subsection 4.1. A non-matching representation of the network is considered by using OPT-XFEM.

$(x, y, z)$ is fixed:

$$h_1(x, y) = \frac{1}{10}\left(-x - \frac{1}{2}\right)\left(8xy\left(x^2 + y^2\right)\operatorname{arctan2}(y, x) + x^3\right),$$
$$h_2(x, z) = \frac{1}{10}\left(-x - \frac{1}{2}\right)x^3 - \frac{4}{5}\pi\left(-x - \frac{1}{2}\right)x^3\,|z|\,,$$
$$h_3(y, z) = (y-1)y(y+1)(z-1)z,$$

in which $\operatorname{arctan2} : \mathbb{R} \times \mathbb{R} \mapsto [-\pi, \pi]$ is the four quadrant inverse tangent function on the 2D reference system $(x_\ell, y_\ell)$ defined as follows:

$$\operatorname{arctan2}(y_\ell, x_\ell) = \begin{cases} \tan^{-1}(y_\ell/x_\ell) & x_\ell > 0 \\ \tan^{-1}(y_\ell/x_\ell) + \pi\operatorname{sign}(y_\ell) & x_\ell < 0 \\ \pi/2\operatorname{sign}(y_\ell) & x_\ell = 0 \end{cases}$$

Function $h$ can be thought as the solution of a Poisson problem with forcing terms $q_i = -\Delta h_i$, on $\Omega_i$, $i = 1, 3$ and Dirichlet boundary conditions $h_D = h(\partial\Omega)$. Figure 5 shows the computational non-conforming triangular mesh, which is representative of the size of elements of the coarsest mesh used for all the considered methods. An example solution is also reported, computed on the same mesh with the optimization method and the XFEM, (OPT-XFEM).

The considered problem, despite being defined on a fairly simple geometry, is interesting as it provides a known analytic solution on a configuration typical of DFN simulations, with multiple fracture intersections, and traces terminating in the inside of the domain. Denote by $\tilde{\Omega}$ the domain without fracture intersections, i.e. $\tilde{\Omega} = \Omega \setminus \left(\bigcup_{k=1}^{3} \Gamma_k\right)$, then the regularity of the solution $h$ is such that $h \in \mathrm{H}^1(\Omega)$, and $h \in \mathrm{H}^2(\tilde{\Omega})$ but $h \notin \mathrm{H}^3(\tilde{\Omega})$. The regularity in $\tilde{\Omega}$ limits the convergence rates for methods conforming to the traces.

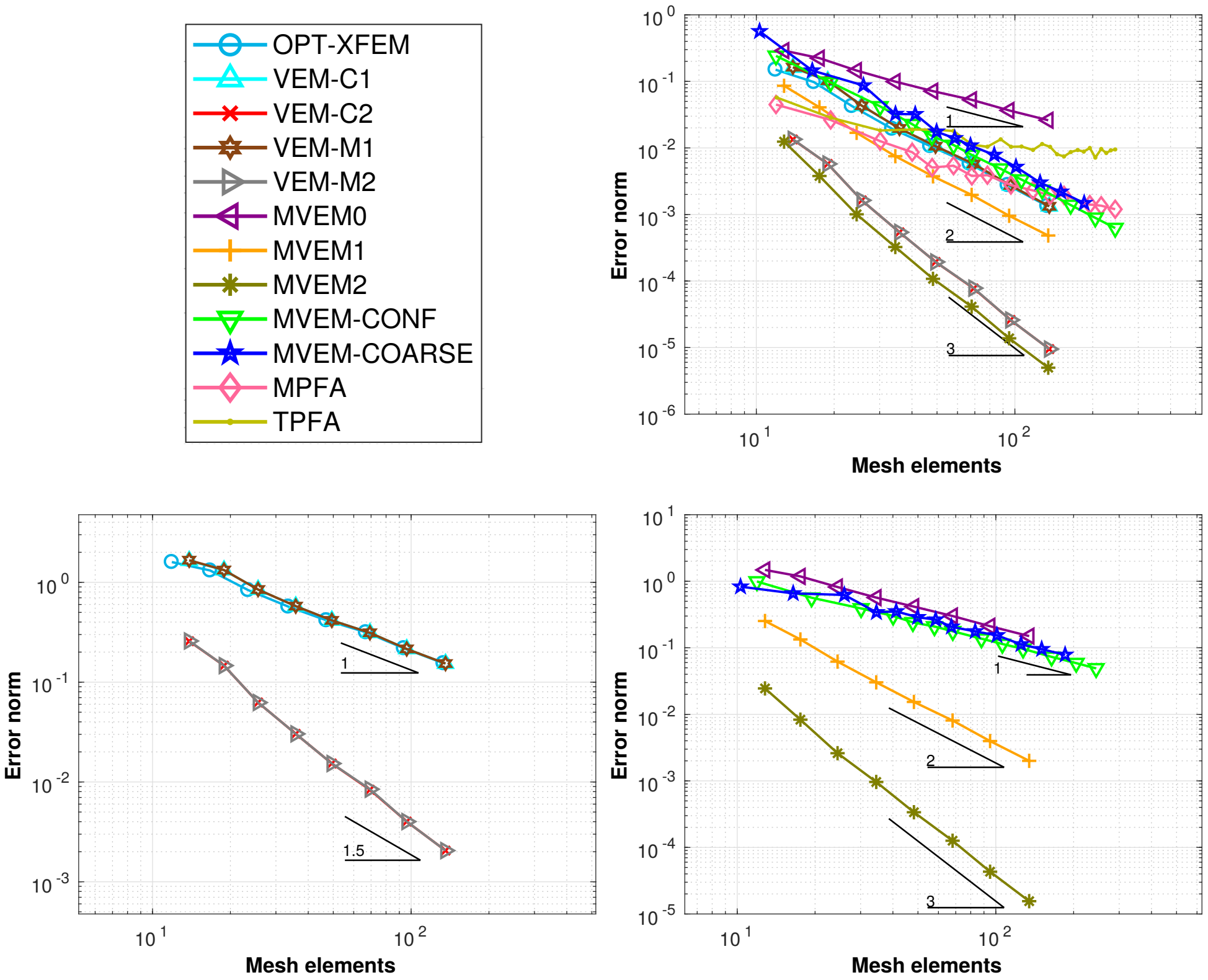

Figure 6: Convergence curves for the setup described in Section 4.1. On the top, the $L^2$ errors in hydraulic head for OPT-XFEM, VEM-C, VEM-M#, MVEM#, MVEM-CONF, MVEM-COARSE, MPFA, and TPFA. On the bottom left, the $H^1$ hydraulic head errors for OPT-XFEM, VEM-C#, and VEM-M#. Finally, on the bottom right the $L^2$ velocity errors for MVEM#, MVEM-CONF, MVEM-COARSE.

To measure the errors, we consider the $\mathrm{L}^2$-norm of the error between the computed hydraulic head in $\Omega$ and the exact solution. Moreover, we compute the $\mathrm{H}^1$ error norms for the methods in primal formulation, and the $\mathrm{L}^2$-norms of the error of both the hydraulic head and of the Darcy velocity with respect to the corresponding exact functions for the methods in mixed formulation.

The convergence results at mesh refinement are reported in Figure 6 (top-right) for the $\mathrm{L}^2$-norm of the error in $h$, in Figure 6 (bottom-left) for the $\mathrm{H}^1$ error norm of $h$ and in Figure 6 (bottom-right) for the error $\mathrm{L}^2$-norm of $-\nabla h$. In all these pictures the error norms are plotted against the number of mesh elements, in order to allow a comparison also among methods in primal and mixed formulation.

Let us start considering the convergence of the hydraulic head, Figure 6 (top-right). We make four observations from the figure. First, we observe that most of the methods convergence as expected. The exceptions are TPFA, which does not converge, and MPFA, which exhibits a somewhat irregular convergence behavior, although the error is consistently decaying. The poor performance of TPFA is expected, as the method is inconsistent for these grids. Moreover, as discussed in Subsection 3.1, the present implementation of MPFA uses a TPFA-type coupling between the fracture planes, and we believe this causes the irregular convergence. It is noteworthy that on the coarsest grids, TPFA and MPFA provides the most accurate solutions of all lowest order methods considered. Second, the methods MVEM-CONF and MVEM-COARSE experience super-convergence, with a convergence order of 2, twice as the expected order for a method with polynomial accuracy order equal to zero. For these two methods, the error has been computed with a zeroth order interpolant of the exact solution, and super-convergence is a known result in this context, see [4]. Third, the optimization-based approach, OPT-XFEM, shows optimal convergence despite the mesh being non-matching. This indicates that the additional basis functions are able to capture the desired irregular behavior across the traces. Finally, the convergence curves of the methods VEM-M and VEM-C of both order 1 and 2 are almost overlapped, since the two matching strategies (Lagrange multipliers and mortaring, re-spectively) have intrinsically a similar nature, with a different definition of the space of the multipliers [9].

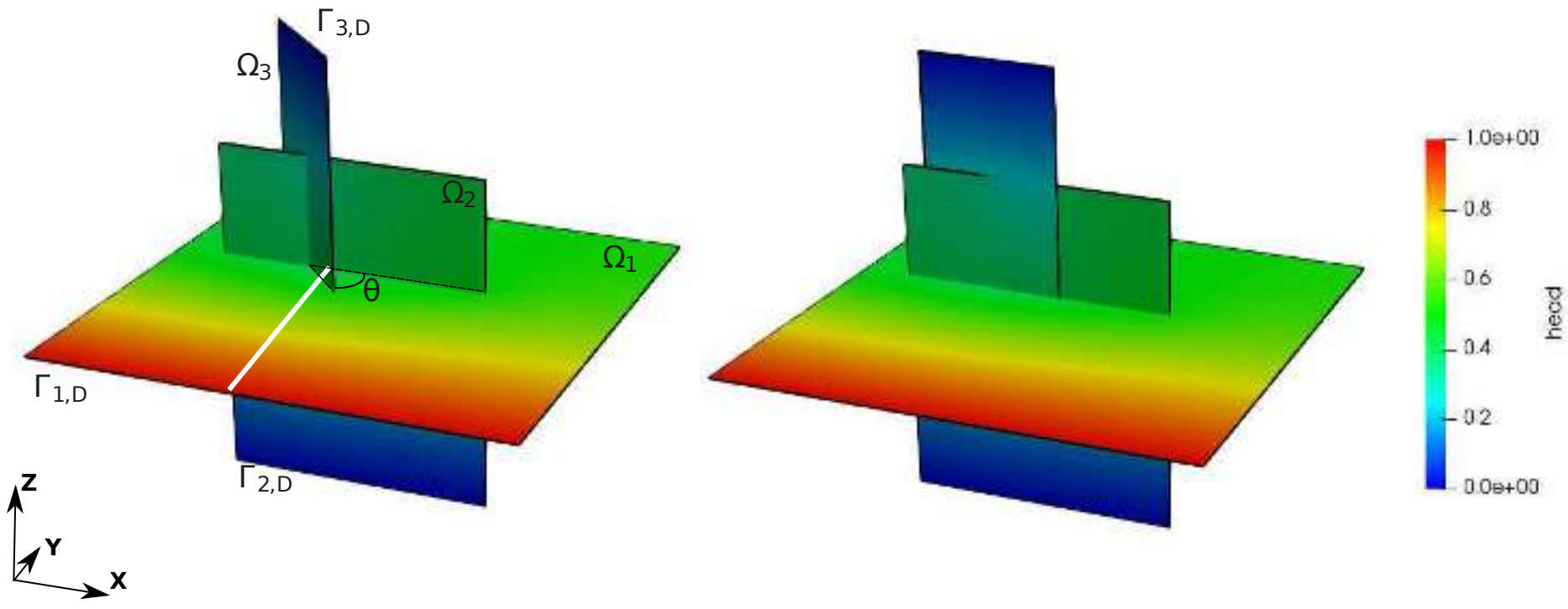

Figure 7: Domain and hydraulic head solution for two examples of Subsection 4.2. In white is indicated the line considered for the plots of Figure 8 and 9.

Next we consider the approximation of the fluxes, computed directly or reconstructed as gradient of the hydraulic head, see Figure 6 (bottom left and right). Again, we see the expected convergence behavior for all methods, with similar observations as for the hydraulic head.

This test confirmed the viability of all methods considered, with the expected caveat for TPFA and partly MPFA. We next go on to test the methods on more challenging geometries.

### *4.2. Vanishing angle between intersecting fractures*

Here, we consider a case where two traces intersect, and study numerical accuracy as the angle between the traces decreases. Such small angles are frequently found in stochastically generated fracture networks, and may also be found in natural networks, depending on the geological setting during network creation. Meshes that conform to this geometry are expected to include a large number of elements, and elongated or badly-shaped elements may be introduced at the intersection between the traces. This may cause problems for methods based on conforming or partly-conforming approaches, while optimization based method should not produce noticeable differences for different angles. We also expect VEM-based approaches to perform well, due to their robustness under rough cell geometries.

The geometry is shown in Figure 7: the DFN is composed of three fractures and three traces, and the two traces on fracture $\Omega_1$ intersect forming an angle $\theta$. Dirichlet boundary condition $h = 1$ is prescribed on the edge marked as $\Gamma_{1,D}$ in Figure 7, homogeneous Dirichlet boundary conditions are set on $\Gamma_{2,D}$ and $\Gamma_{3,D}$. All other fracture edges are assigned homogeneous Neumann conditions, and no forcing terms are present. We consider a series of angles $\theta$, ranging between $\pi/2$ and $\pi/565$ with 20 equally spaced steps. A computed numerical solution is also shown in Figure 7 for two values of the angle $\theta$: $\theta = \pi/3$ on the left and $\theta = \pi/180$ on the right.

We solve the above sequence of problems with the proposed numerical methods. Two sets of grids are considered, where the finer grid has about twice as many cells as the coarser. For grids that are conforming to the traces, the number of cells will blow up as the angle between the traces closes. Non-matching grids will show no such behavior.

To measure the accuracy of the methods, we consider the outflow from the network, as a proxy for the upscaled permeability, and also the trace of the solution on a segment $\gamma$ placed on fracture $\Omega_1$ defined as

$$\gamma = \{(x, y, z) : x = 0.35, y \in (0, 1), z = 0\},$$

see Figure 7 for the reference coordinate system. Results are shown in Figures 8 on the coarser mesh and in Figure 9 for a finer mesh. The top frames of both figures show the line plot of computed hydraulic head solution for the two extreme values of the angle $\theta$, with $\theta = \frac{\pi}{2}$ on the top-left and $\theta = \frac{\pi}{565}$ on the top-right. All solutions are in good agreement for both grid levels. There is larger difference in the total outflow from the network: considering the method based on the VEM on conforming triangular mesh (label MVEM-CONF) as the reference, the methods either slightly overestimate (OPT-FEM, OPT-XFEM, VEM-C#, VEM-M#) or slightly underestimate (TPFA, MPFA and MVEM-COARSE) the flux. The larger differences are produced by the OPT-FEM and TPFA methods which however differ from the reference solution of less than 5% on the coarse mesh and less than 3% on the finer mesh. The differences among all the methods decrease with grid refinement.

Perhaps the most interesting aspect of this test is the number of cells needed to mesh the network according to method requirements, shown in Figure 8 and Figure 9, both lower left. The cell count for the matching meshes increases drastically as the angle closes, while for the non-matching cases, the cell count is more or less constant for all angles. The matching grids were created by standard meshing software, like Gmsh [34]. To mimic simulations in larger scale networks, where manual tuning of meshes in all difficult regions may not be feasible, no attempts were made to guide the meshing software apart from setting general mesh size requirements. The high cell count is thus an unavoidable effect of conforming meshing in the presence of almost parallel traces. The method MVEM-COARSE post-processes the standard conforming mesh, gluing the original triangles into conforming polygons, and is thus able to partly mitigate the increase of mesh cell number for the smaller values of $\theta$.

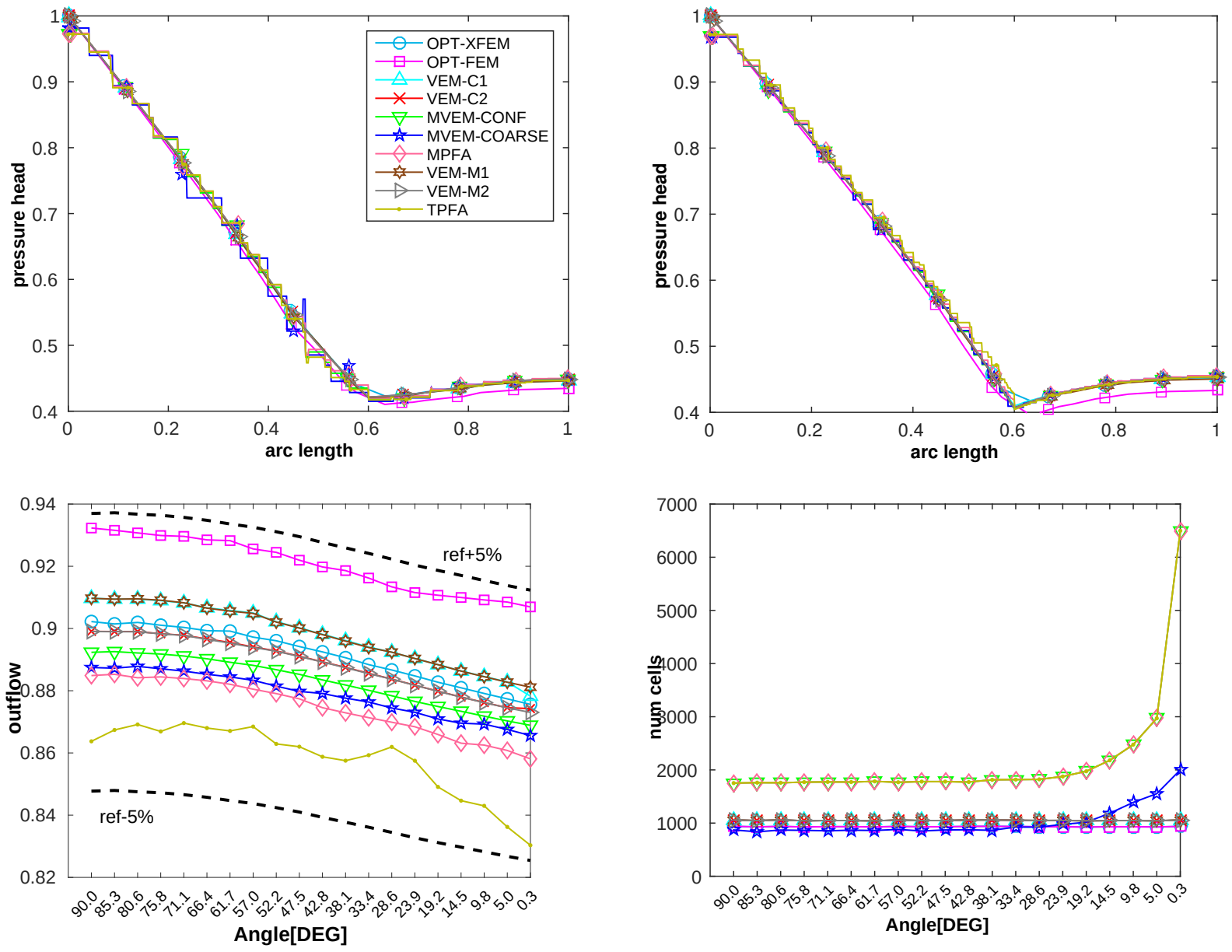

Figure 8: Coarse mesh results for example in Subsection 4.2. On top the hydraulic head over $\gamma$, on the bottom left the outflow and the right the number of cells.

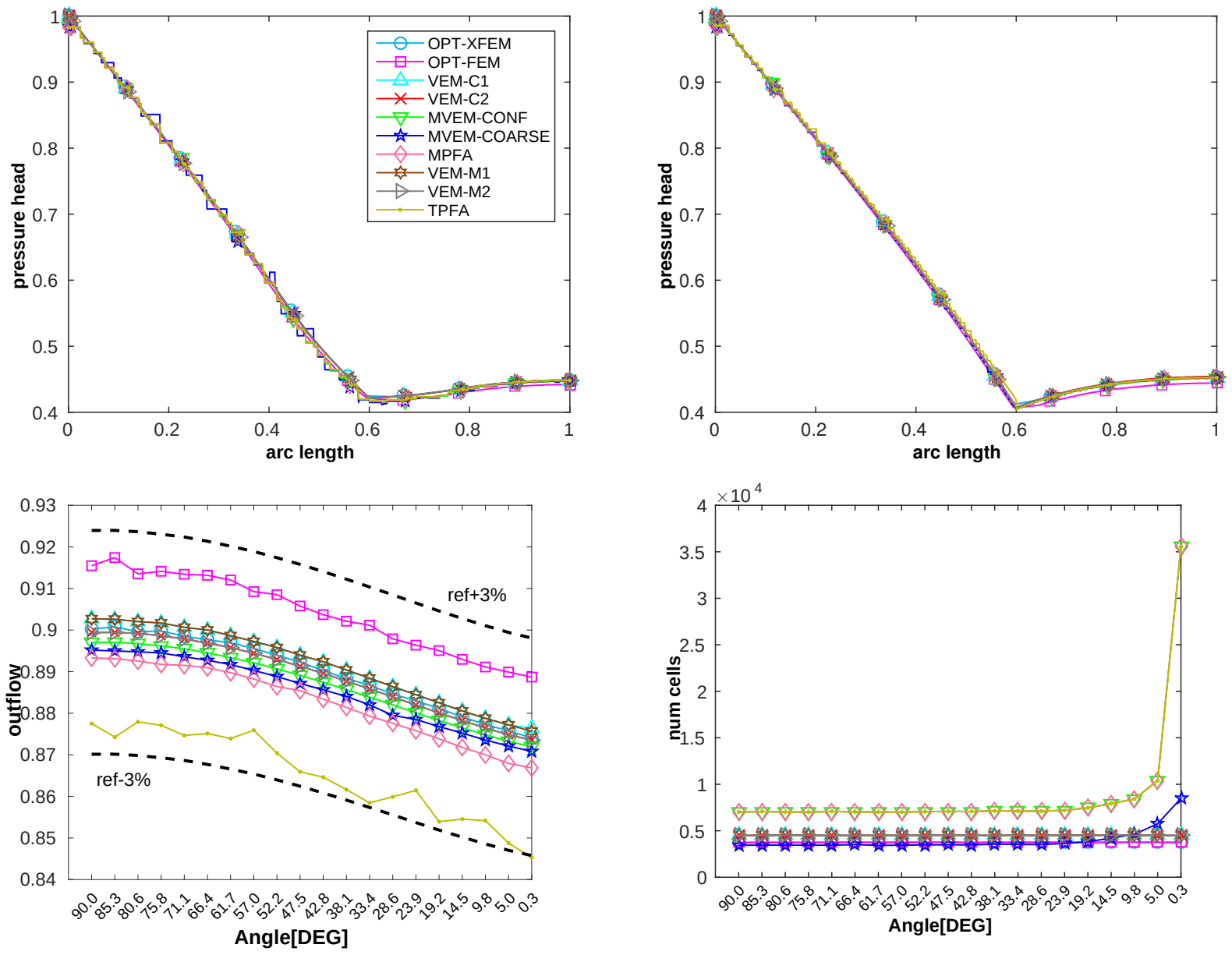

Figure 9: Fine mesh results for example in Subsection 4.2. On top the hydraulic head over $\gamma$, on the bottom left the outflow and the right the number of cells.

### *4.3. Vanishing trace between intersecting fractures*

In this part, we consider the impact of a fracture sliding out a network. The setup generates small fracture traces, and the test aims to analyze the accuracy and the stability of the presented numerical approaches in this setting. Small traces may impose challenges associated to the loss of regularity of the solution. We propose two tests, one with a simple geometry in Subsection 4.3.1 and one with a more complex geometry in Subsection 4.3.2, and modify the geometry in both cases so that the trace length decreases. In the first test, the vanishing trace is a chock point, while in the second the flow can reroute to other parts of the network. Both effects are important for practical networks. In both cases we present plots over a line for the hydraulic head for several configurations as well as the net flux through the network, comparing the results of the numerical schemes.

#### *4.3.1. Flow break due to vanishing trace*

For this first case we consider a network composed by three fractures, where their relative position changes over 21 different configuration. Referring to Figure 10 the fracture on the left, called $\Omega_l$, is fixed while the other two fractures, the middle $\Omega_m$ and the right $\Omega_r$, are moving on the right for each configuration: starting from configuration id 1 (leftmost in Figure 10) the trace between fracture $\Omega_m$ and $\Omega_l$ ranges from a length of 0.6 to a length of 0.01 at configuration id 21 (rightmost in Figure 10). For simplicity we assume unitary permeability in all the fractures, unitary hydraulic head for the bottom boundary of $\Omega_l$ and zero hydraulic head for the bottom boundary of $\Omega_r$. The other boundaries are impervious. For all the methods the number of cells for the finest mesh is approximately 1100 per fracture.

The computed hydraulic head in this DFN is reported in Figure 10 for configuration ids 1 (left), 10 (center) and 21 (right) on the finest mesh, giving an idea of its behavior at changing the geometry. In the images, the gradient in hydraulic head on $\Omega_l$ becomes steeper as the trace between $\Omega_l$ and $\Omega_m$ becomes smaller.

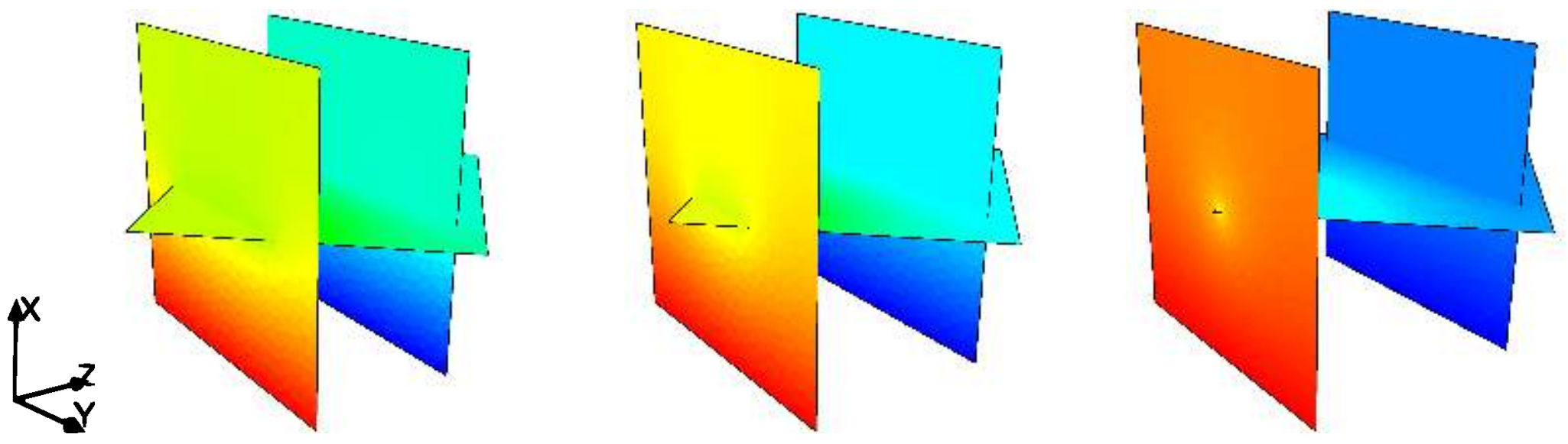

Figure 10: Representation of three geometrical configurations of the network for the example described in part 4.3.1. The pressure solution ranges from 0 to 1.

To compare the behavior of the numerical schemes presented in the previous sections, we consider three plots over line for the hydraulic head at the line $\gamma$ on $\Omega_l$, defined as:

$$\gamma = \{(x, y, z) \in \mathbb{R}^3 : x = 0.5, y = 0.5, z \in (z_{\min}, z_{\max})\},$$

where $z_{\min}$ and $z_{\max}$ are the minimum and maximum of the $z$-coordinate for the network, respectively. Note that $\gamma$ changes position at each configuration. The line $\gamma$ lies on the middle line of $\Omega_m$, thus we consider the arch length as abscissa for the plots over $\gamma$.

Figure 11 (top left and right and bottom left) shows the comparison between the numerical schemes of the hydraulic head over $\gamma$ for the geometry with id 1, 10, and 21, respectively, on the finest mesh. The solutions in the graphs are in good agreement irrespective on the coupling strategies between fractures adopted by the various methods. The only noticeable difference is with OPT-FEM for the last configuration, which produces a slightly higher value of the hydraulic head at the left of the arch. The stair shape behavior of some methods is due to the piece-wise approximation of the hydraulic head in each element, which causes graphical oscillations. Finally, we compare the outflow from the network for the different geometries among the numerical schemes, see Figure 11 (bottom right). Taking as reference the solution obtained with the method MVEM-CONF on this fine mesh, we observe that also in this case the solutions are in good agreement when the length of the vanishing trace is greater than $\sim 0.07$. For smaller values some of the methods give differences slightly higher than 10% of the reference. In particular the larger differences are given by methods OPT-FEM which overestimates the flux and MVEM-COARSE, which instead underestimates it.

In Figure 12 we study the same problem but on two different coarser meshes, having about 150 and 30 cells per fracture, respectively, and we compare both the hydraulic head at the geometrical configuration id 21 and the flux for all the simulations. For the hydraulic head (Figure 12 left), a poorer resolution is clearly visible, but mainly for the methods providing a piece-wise approximation of the solution. Considering the fluxes (Figure 12 right), on the intermediate mesh many of the methods still remain close to the reference, the largest differences produced again by the methods OPT-FEM and MVEM-COARSE when trace length becomes smaller than about 0.15 and by VEM-C and VEM-M of order 1 for trace lengths smaller than $\sim 0.07$. On the coarsest mesh, all the considered methods give differences higher than 10% of the reference for some values of the length of the vanishing trace. The best results are produced by methods TPFA, MPFA, MVEM-CONF, OPT-XFEM and VEM-C and VEM-M of order 2. The method MVEM-COARSE, even if produces results with differences slightly larger than 10% for intermediate values of trace length is however capable of catching the overall behavior of the solution. Both MPFA and MVEM with and without coarsening, shows somewhat irregular behavior for coarser grids. It is to remark that this test is extremely challenging for polygonal/non-conforming/non-matching methods. In fact, conforming meshes are necessarily smaller close to the small trace, where the solution has the steepest gradients, whereas non-matching meshes (as the one for OPT-XFEM and OPT-FEM) are arbitrarily placed with respect to the trace. A similar argument applies to the polygonal meshes used for VEM-M and VEM-C methods and for the method MVEM-COARSE, whose mesh is a non-guided agglomeration of a former triangular conforming mesh, see Figure 13. These methods are in fact designed to be extremely low computationally demanding and robust to hard to mesh geometries. Methods as OPT-XFEM, VEM-C and VEM-M of high order partly compensate the non conformity of the mesh with enlarged function spaces.

*4.3.2. Flow redistribution due to vanishing traces*

This example considers a more complex geometry composed by 10 fractures, where a fracture $\Omega_s$ shrinks over 44 different configurations: the length of fracture $\Omega_s$ ranges from 2 (configuration id 1) to 0.26 (configuration id 44). As this fracture gets smaller, less connections are present in the system. Figure 14 represents some configurations, $\Omega_s$ being the horizontal middle fracture. We assume unitary permeability in all the fractures. Boundary conditions are specified as hydraulic head of one and zero on the bottom boundary of the vertical fractures on the left and right, respectively. The other boundaries are impervious. The mesh for this example is obtained setting a minimum number of $\sim 500$ cells per fracture.

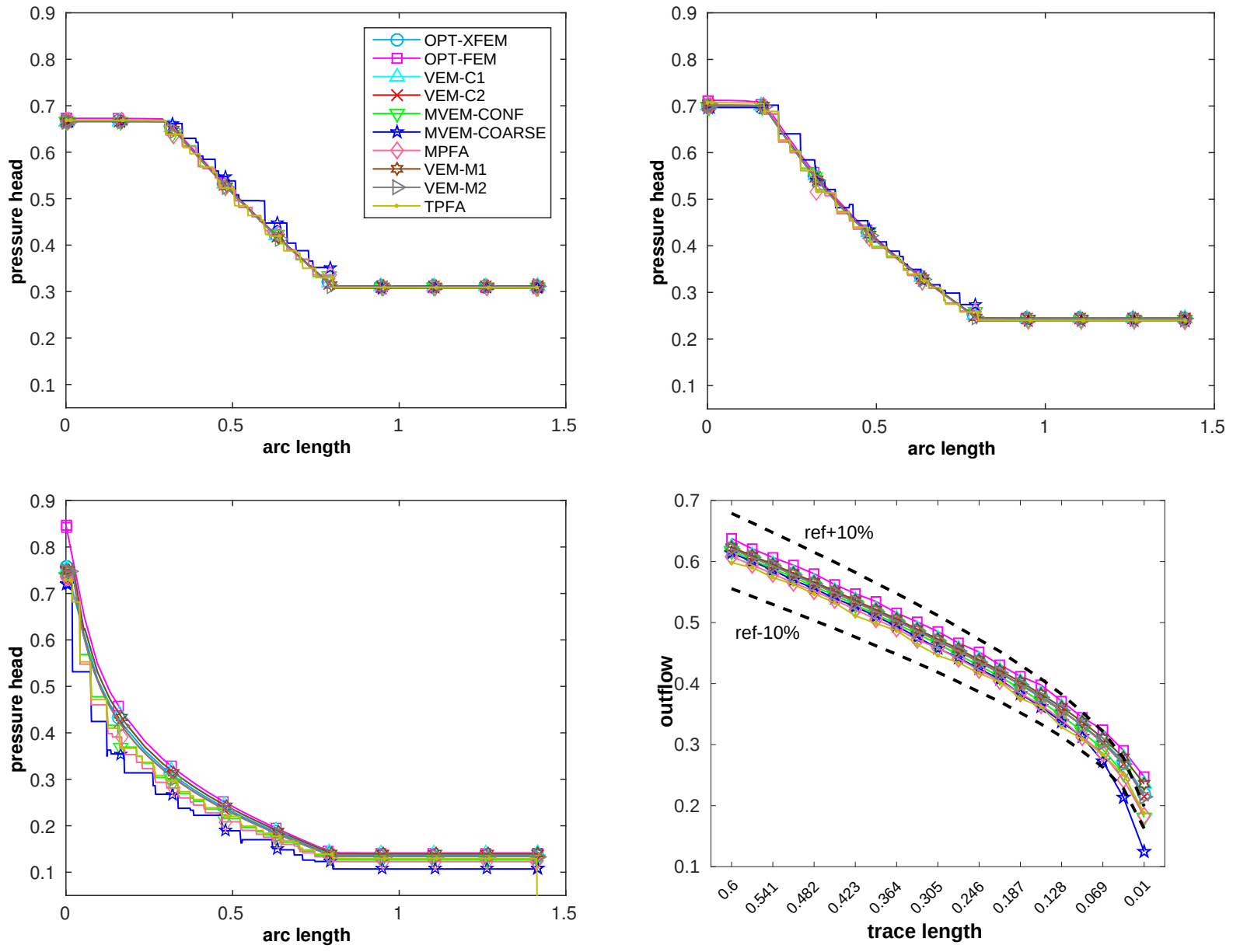

Figure 11: The two plots on top and the one on bottom left represent the hydraulic head over line $\gamma$ for the example in part 4.3.1. The plot on the bottom right depicts the outflow in function of the geometry.

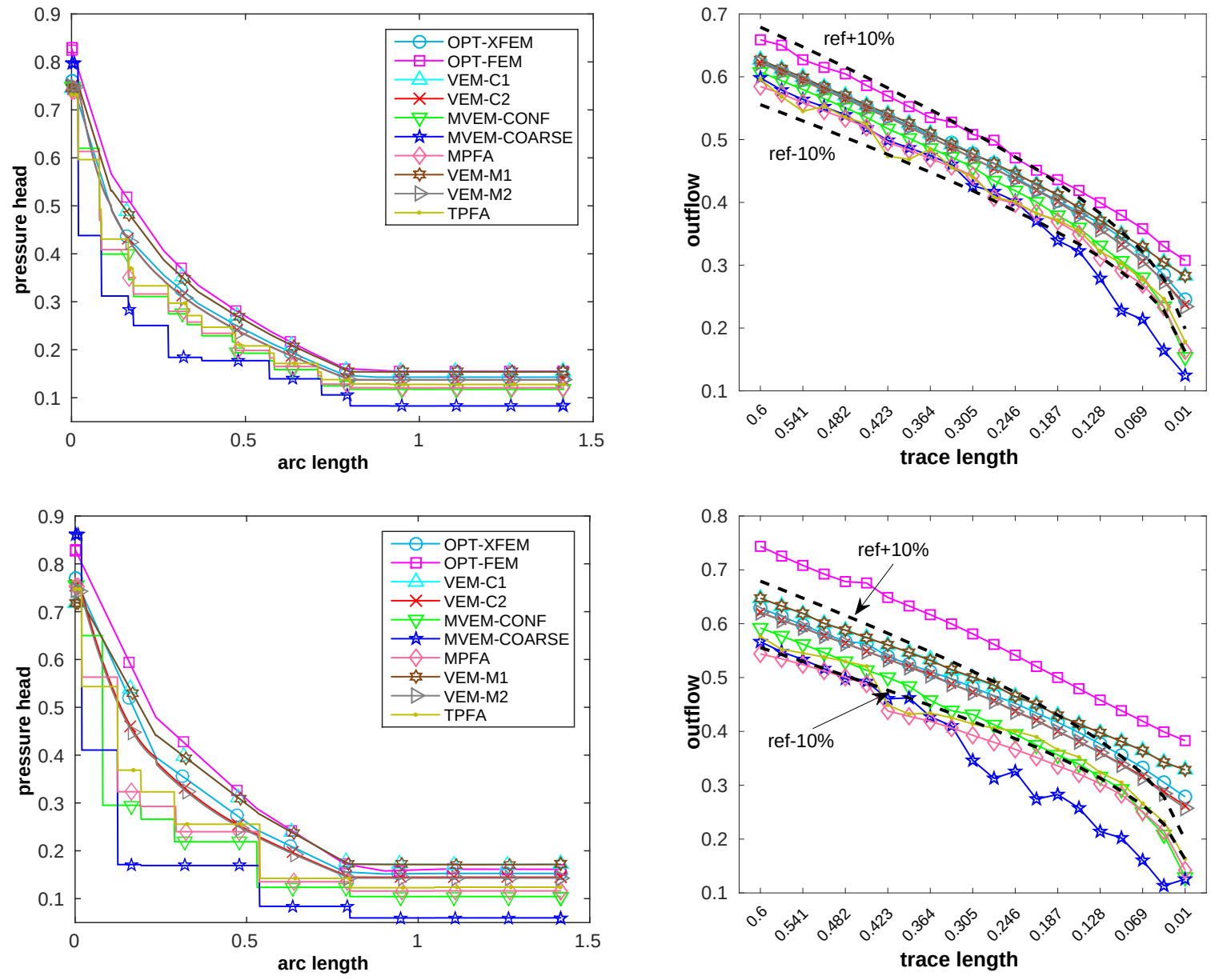

Figure 12: Pressure over $\gamma$ and flux comparison for the example in part 4.3.1. Coarse meshes with $\sim$150 and $\sim$30 elements per fracture on the top and bottom, respectively. Left column: hydraulic head along $\gamma$; right: outflow from the network.

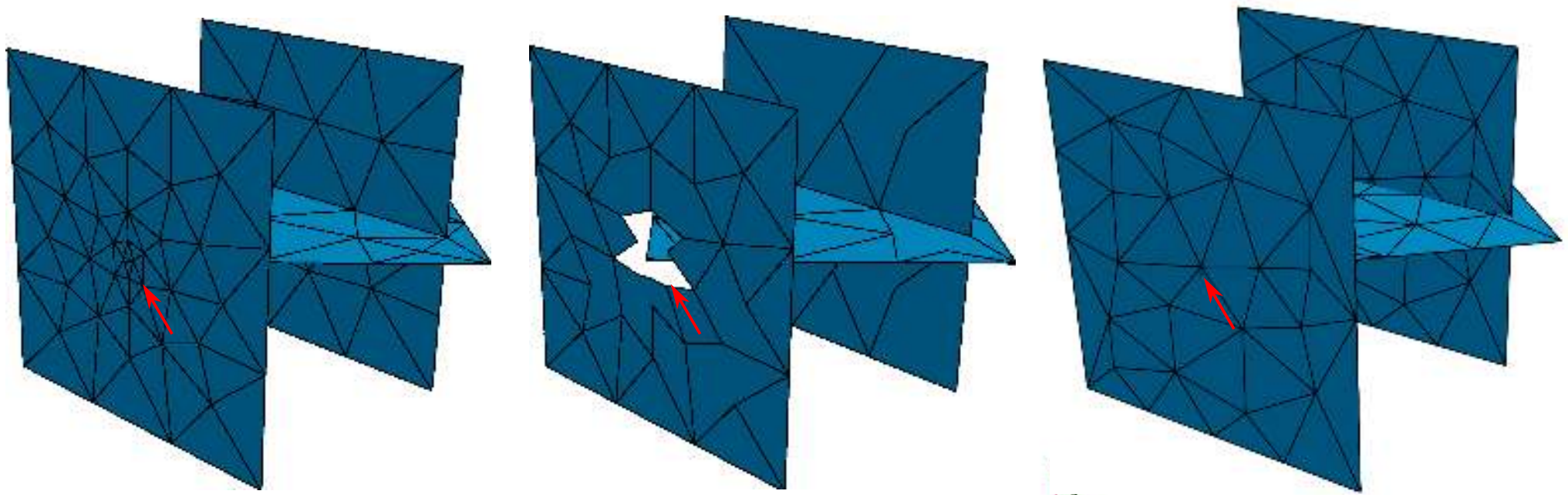

Figure 13: Coarse mesh ($\sim 30$ elements per fracture) for the example in part 4.3.1: conforming (left), coarsened (center) non-matching (left). For the coarsened mesh, the missing element is due to the visualization, the polytopal cell has internal edges.

The previous case considered a network that in the limit case becomes disconnected. In this case the problem is designed in a way that, as $\Omega_s$ is reduced the flux can be redistributed on other fractures. The test thus allows us to study methods behavior in a limiting case, but now allowing for flow bypass of the critical points.

The solution is represented in Figure 14 for some configurations. To compare the behavior of the numerical schemes we consider, also in this case, a plot over line for the hydraulic head for each considered geometry. The line $\gamma$ is now defined as

$$\gamma = \{(x, y, z) \in \mathbb{R}^3 : x = 1.5, y \in (0, 2), z = 0.5\},$$

which belongs to the vertical fracture on the right of Figure 14. Figure 15 represents the comparison among the numerical schemes of the hydraulic head over $\gamma$ for configuration id 1 (prior to sliding), on the top-left, and configuration id 44 (when $\Omega_s$ is almost disconnected), on the top-right. The hydraulic heads shown in each plot are in good agreement between the methods.

In Figure 15 (bottom right) we compare the outflow from the system at different geometrical configurations. Most of the methods behave similarly with a jump in the flux every time $\Omega_s$ detaches from one of the vertical fractures. The main differences are represented by OPT-FEM and MPFA. In the former the flux is over estimated while in the latter is under estimated, but however confined in a difference of 5% with respect to the solution obtained with the method MVEM-CONF, considered as a reference, and the general behavior of the solution is well reproduced. As already discussed the method OPT-FEM is built on non matching meshes and no additional basis functions are introduced to compensate the non conformity, as to keep the computational cost of the method low. In the case of the multi-point flux approximation the quality of the solution can be affected by the implementation of the coupling conditions between fractures with the two-point flux approximation, which in this more complex case have a more relevant impact. Compared to the simpler case, the differences in fluxes between the methods are systematic as the geometry is modified, and probably indicate general approximation properties, rather than specifics related to this setup. Thus, as expected, small traces are less of a problem when the flow can bypass the choke point, and also polygonal/non-conforming/non-matching methods can produce higher quality results.

### *4.4. Upscaling for extruded real outcrop*

Our final test considers a network generated from a real outcrop, located in Western Norway (GPS coordinates 60°20′19.8″N and 4°55′54.7″E), see Figure 16. From the interpreted fractures, we pick out a region containing 66 fracture lines, and extrude them to 3D fractures in a way that correlates fracture radius with the observed line length, and also preserves abutting relations between the fractures (fractures that meet in a T-intersection in the outcrop obey the same relation in 3D). On the resulting network, we impose a bounding box so that a few fractures hit each face of the box. We then calculate the bulk flow response by imposing a unit drop in hydraulic head in one coordinate direction at the time, assigning no-flow conditions on the other boundaries and having unitary permeability.

Contrary to the previous examples and due to the number of fractures, we consider two ways for constructing conforming grids for MVEM-CONF, MVEM-COARSE, MPFA, and TPFA. The first considers a conforming discretization on each trace as done previously, in the second each fracture is meshed independently and then conformity is restored splitting the edges of the elements on the two fractures meeting at each trace. We expect less cells by considering the latter case, thanks to the independent meshing process. In the results, we indicate it by full purple markers of the same shape of the related method.

The outflow from the fracture network in each direction, plotted as functions of the number of cells, is shown in Figure 17 for the three principal directions. The values converge as the grids are refined, and there is a spread of about $\sim 10\%$ in all the results, between the estimated outflows for the finest grids. The method VEM-M of order 1 shows the larger discrepancies with respect to the average values obtained by the other methods. This can be explained by the fact that no mesh adjustments where made in order to meet quasi-uniformity requirement on the mesh, necessary to ensure the convergence of the Mortar multipliers to the flux [7].

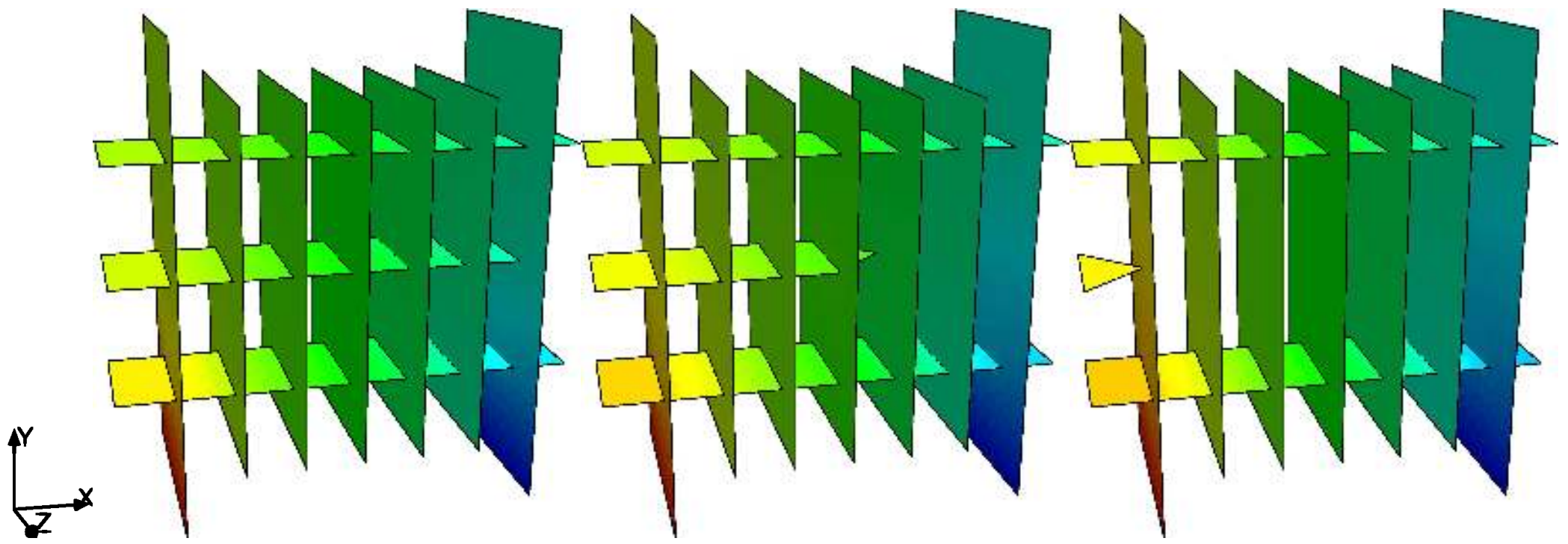

Figure 14: Representation of four geometrical configurations of the network for the example described in Subsection 4.3.2. The computed hydraulic heads range from 0 to 1.

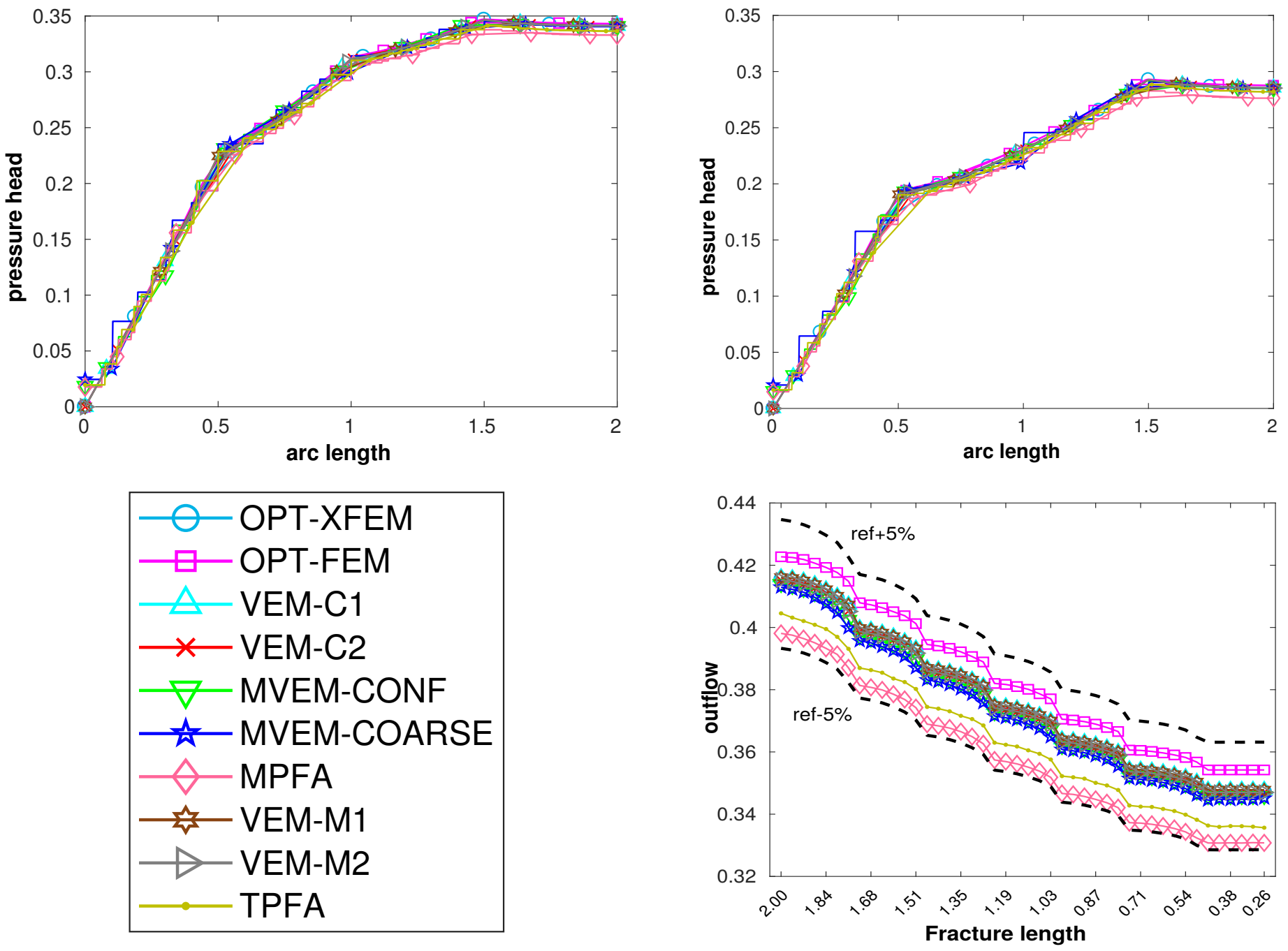

Figure 15: On the top, plot over line of the hydraulic head for the first (left) and last (right) simulation form the example of Subsection 4.3.2. On the bottom, outflow with respect to the fracture lengh.

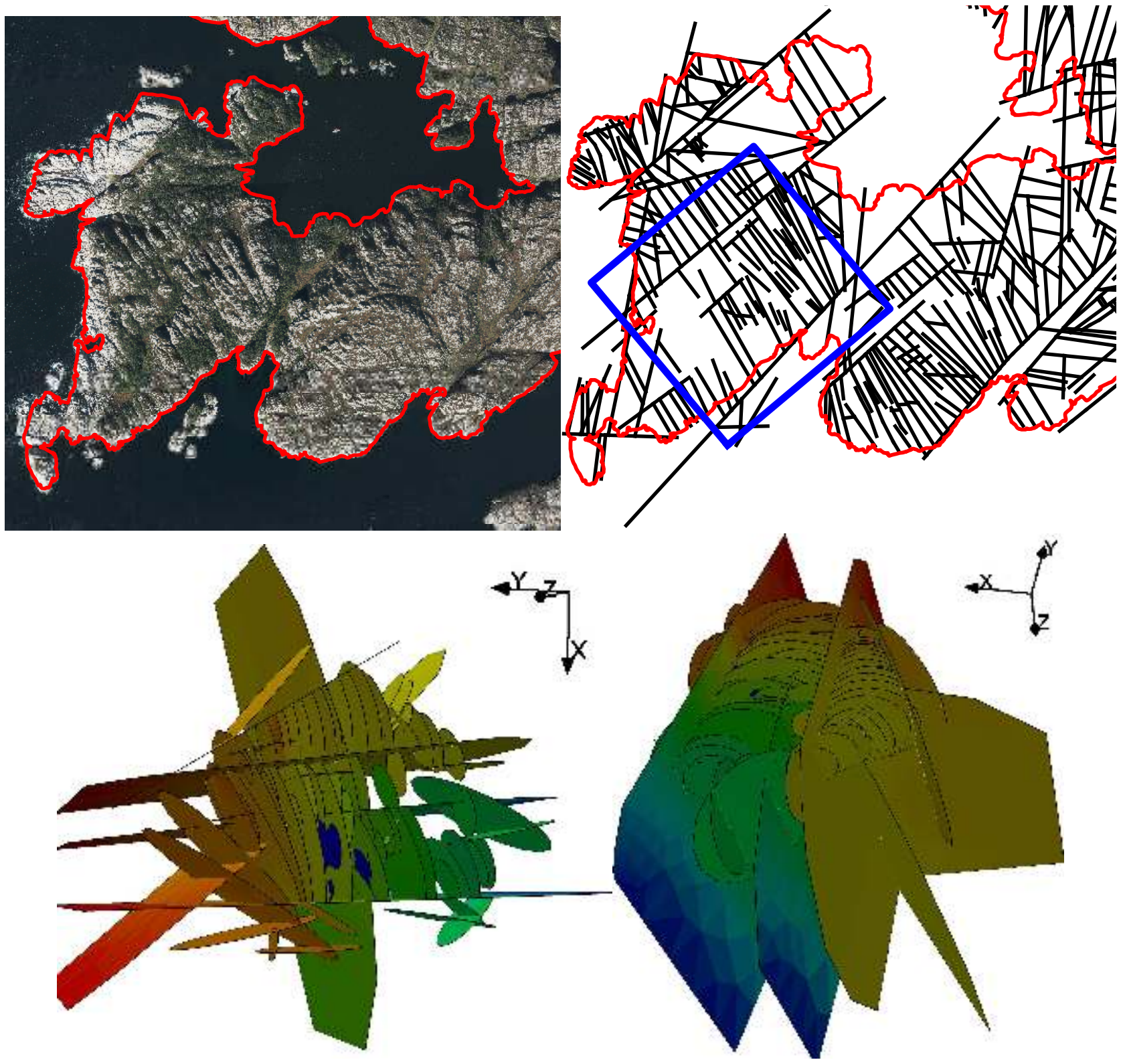

Figure 16: On the top, the map and the interpreted outcrop with the region of interest marked with a blue square. On the bottom, 3D views of the extruded fractures colored by the computed hydraulic head for one case. Some of the fractures are blue since they are disconnected from the main network.

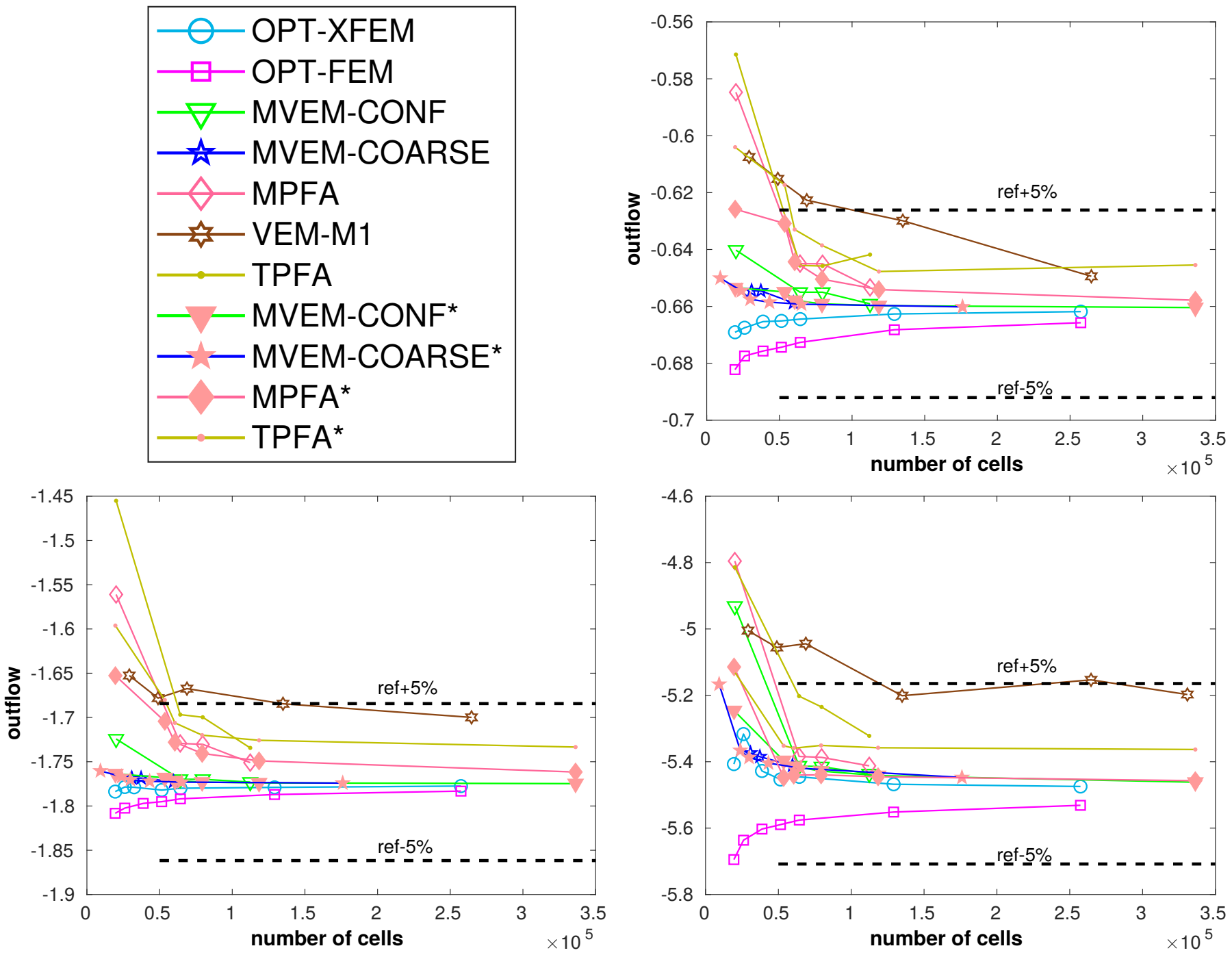

Figure 17: Flux computation for the three directions. On the top along the $x$ axis, on the bottom left along the $y$ axis, and on the bottom right along the $z$ axis.

Let us consider as reference the solution obtained with the MVEM-CONF on the finest mesh. In practical simulations, there is often no upper limit on how many fractures should ideally be included, or correspondingly how fast each simulation should run. To that end, a more interesting test is to consider how the results diverges as the number of cells decreases. In this interpretation, the most stable methods are OPT-XFEM, MVEM-CONF* and MVEM-COARSE. The other methods diverge significantly for coarser grids on all flow directions. Even for the best methods, we observe some irregularities on the coarsest grids for a drop in hydraulic head in the $z$-direction. Naturally, there is a limit to how cheap the simulations can be, independent of the discretization applied, however, characterizing this limit is non-trivial.

## 5. Conclusions

Discretization methods for flow in DFNs have received considerable attention in recent years. The proposed numerical methods differ in their approach to meshing of the intersections between fractures, the coupling of flow between fractures, and the discretization schemes applied. This work has reviewed several approaches to all of these ingredients, and applied the methods to a range of test problems. Overall, the methods were in relative good agreement. Differences in the upscaled network permeability consistently showed larger errors than localized differences in hydraulic head. The proposed numerical tests revealed that methods based on polygonal/non-conforming/non-matching meshes can produce results comparable to those of methods based on conforming meshes, with less limitations related to mesh generation. There is a clear trade-off between computational cost and accuracy: when non conforming meshes are used, increasing the discrete functional spaces improves the performances. In particular, good performances were shown by an approach based on the optimization and the XFEM, which allows for non-matching meshes and couples the solution between fractures via an optimization approach, and a method based on the mixed VEM, that can deal with partly conforming meshes and general cell shapes. In both synthetic and more realistic problems, these methods consistently showed good accuracy in particular for coarser meshes. Polygonal and non-matching mesh based approaches thus seem to be the viable option for flow in large-scale networks, where computational cost is the limiting factor for the simulations.

## Acknowledgment

The first and second authors acknowledge financial support for the ANIGMA project from the Research Council of Norway (project no. 244129/E20) through the ENERGIX program. The third author acknowledges financial support from GNCS - INdAM and from Politecnico di Torino through project "Starting grant RTD". The authors warmly thank Luisa F. Zuluaga for constructing and providing the real fracture network for the example in Subsection 4.4. Finally, the authors warmly thanks Ivar Stefansson for many fruitful discussions related to the development of this work.